\documentclass[a4paper,11pt]{article}

\usepackage[margin=3cm,top=3cm,footskip=15mm]{geometry}

\usepackage{amsmath, amssymb, amsthm, mathrsfs}
\usepackage[dvipsnames]{xcolor}

\usepackage[colorlinks=true,linkcolor=red!70!black,citecolor=MidnightBlue]{hyperref}

\usepackage[utf8]{inputenc} 
\usepackage[T1]{fontenc}
\usepackage{lmodern}

\numberwithin{equation}{section}
\newcommand{\N}{\mathbf N}
\newcommand{\R}{\mathbf R}

\def\AA{{\mathcal A}}
\def\BB{{\mathcal B}}

\def\HH{{\mathcal H}}

\def\MM{{\mathcal M}}

\def\OO{{\mathcal O}}

\def\RR{{\mathcal R}}

\def\VV{{\mathcal V}}

\def\VV{{\mathcal V}}

\def\CCC{{\mathscr C}}

\def\LLL{{\mathscr L}}

\def\RRR{{\mathscr R}}

\def\le{\leqslant}
\def\ge{\geqslant}

\DeclareMathOperator{\Div}{div}

\def\eps{{\varepsilon}}
\def\wto{{\,\rightharpoonup\,}}

\newcommand{\la}{\left\langle}
\newcommand{\ra}{\right\rangle}

\newcommand{\Nt}{|\hskip-0.04cm|\hskip-0.04cm|}

\renewcommand{\d}{\mathrm{d}}
\newcommand{\dv}{\mathrm{d} v}
\newcommand{\dt}{\mathrm{d} t}

\newcommand{\dx}{\mathrm{d} x}

\newtheorem{theo}{Theorem}[section]
\newtheorem{prop}[theo]{Proposition}
\newtheorem{lem}[theo]{Lemma}

\theoremstyle{definition}
\newtheorem{rem}[theo]{Remark}

\newcommand{\be}{\begin{equation}}
\newcommand{\ee}{\end{equation}}
\newcommand{\ba}{\begin{aligned}}
\newcommand{\ea}{\end{aligned}}
\newcommand{\beqn}{\begin{equation}}
\newcommand{\eeqn}{\end{equation}}
\newcommand{\bear}{\begin{eqnarray}}
\newcommand{\eear}{\end{eqnarray}}
\newcommand{\bean}{\begin{eqnarray*}}
\newcommand{\eean}{\end{eqnarray*}}

\allowdisplaybreaks

\title{Hypocoercivity for kinetic linear equations in bounded domains with general Maxwell boundary condition}
\author{Armand Bernou, Kleber Carrapatoso, St\'ephane Mischler and Isabelle Tristani}

\def\signab{(A.~Bernou) Sorbonne Universit\'e and Universit\'e de Paris, CNRS, Laboratoire Jacques-Louis Lions (LJLL), F-75005 Paris, France.

\textit{Email address}: \texttt{armand.bernou@sorbonne-universite.fr} \bigskip}

\def\signkc{(K.~Carrapatoso) CMLS, \'Ecole polytechnique, Institut Polytechnique de Paris, 91128 Palaiseau Cedex, France.

\textit{Email address}: \texttt{kleber.carrapatoso@polytechnique.edu} \bigskip}

\def\signsm{(S.~Mischler) CEREMADE (CNRS UMR 7534), PSL University, Université Paris-Dauphine, Place de Lattre de Tassigny, 75775 Paris 16, France. 

\textit{Email address}: \texttt{mischler@ceremade.dauphine.fr} \bigskip}

\def\signit{(I.~Tristani) D\'epartement de Math\'ematiques et Applications, \'Ecole Normale Sup\'erieure, CNRS, PSL University, 75005 Paris, France. 

\textit{Email address}: \texttt{isabelle.tristani@ens.fr}}

\begin{document}

\maketitle
\begin{abstract}
We establish the convergence to the equilibrium for various linear collisional kinetic equations (including linearized Boltzmann and Landau equations) with physical local conservation laws in bounded domains with general Maxwell boundary condition. 
Our proof consists in establishing an hypocoercivity result for the associated operator, in other words, we exhibit a convenient Hilbert norm for which the  associated operator is coercive in the orthogonal of the global conservation laws. Our approach allows us to treat general domains with all type of boundary conditions in a unified framework. In particular, our result includes the case of vanishing accommodation coefficient and thus the specific case of the specular reflection boundary condition.   
\end{abstract}

\vspace{-.3cm}

\tableofcontents

\section{Introduction}
\label{sec:Intro}

\subsection{The problem} \label{subsec:equation}
In this paper, we study a linear collisional kinetic equation in a bounded domain with general Maxwell boundary condition. More precisely, we consider a smooth enough bounded domain $\Omega \subseteq \R^d$, $d \ge 2$, we denote by $\OO := \Omega \times \R^d$ the interior set of phase space and $\Sigma := \partial\Omega \times \R^d$ the boundary set of phase space. 
For a (variation of a) density function  $f=f(t,x,v)$, $t \ge 0$, $x \in \Omega$, $v \in \R^d$, we then look at the following equation 
\bear\label{eq:dtf=Lf}
\partial_t f &=& \LLL f := - v \cdot \nabla_x f +  \CCC f \quad\hbox{in}\quad (0,\infty) \times \OO, 
\\ \label{eq:BdyCond}
\gamma_{\!-} f &=& \RRR \gamma_{\!+}  f \quad\hbox{on}\quad (0,\infty) \times \Sigma, 
\eear
where $\gamma_{\!\pm} f$ denote the trace of $f$ at the boundary set and where $\CCC$ and $\RRR$ stand for two linear collisional operators that we describe below. Our goal is to investigate the long-time behavior of solutions to this linear equation. In order to do so, we will prove an hypocercivity result using a general and robust approach inspired by previous works on~$L^2$-hypocoercivity. 

\smallskip
\noindent{\it Motivation.} We first briefly explain the motivation to study this problem. We consider a system of particles confined in~$\Omega$ whose state is described by the  variations of the density of particles $F = F(t,x,{v}) \ge 0$ which at time $t \ge 0$ and at position~$x \in \Omega$, move with velocity ${v} \in \R^d$. We suppose that collisions between particles are for instance described by the Boltzmann or the Landau bilinear collision operator. It leads us to consider the following equation:
\bear \label{eq:nonlinear}
\partial_t F &=& - v \cdot \nabla_x F + Q(F,F) \quad\hbox{in}\quad (0,\infty) \times \OO \\ \label{eq:Bdynonlinear}
\gamma_{\!-} F &=& \RRR \gamma_{\!+}  F \quad\hbox{on}\quad (0,\infty) \times \Sigma, 
\eear
where $Q$ is for instance the Boltzmann or the Landau collision operator. The standard (normalized and centered) Maxwellian 
\beqn\label{eq:standardMaxw}
\mu=\mu(v) := (2\pi)^{-d/2}e^{-|v|^2/2}
\eeqn
is a global equilibrium of this equation. 
In order to study this type of problem in a close-to-equilibrium regime, we write the distribution $F$ as the following perturbation of the global equilibrium $\mu$: $F=\mu+ f$. If $F$ solves~\eqref{eq:nonlinear}-\eqref{eq:Bdynonlinear}, then the linearized equation (throwing away the quadratic term) satisfied by $f$ is nothing but \eqref{eq:dtf=Lf}-\eqref{eq:BdyCond} with 
$$
\CCC f := Q(\mu,f) + Q(f,\mu).
$$ 
The assumptions (A1)-(A2)-(A3) made below on the collisional operator~$\CCC$ are met by the linearized Boltzmann and Landau equations for the so-called hard potentials (and thus including the Boltzmann hard spheres case). It is worth noting that by a straightforward adaptation of our method, we can also treat linear operators preserving only mass such that the Fokker-Planck operator or the relaxation operator. We believe that our analysis is also new in this setting.
In Section~\ref{sec:weak}, we present more general assumptions that allow us to deal with linearized Boltzmann and Landau operators corresponding to softer potentials.

\smallskip
\noindent{\it The boundary condition.} Let us now describe the boundary condition \eqref{eq:BdyCond}. For that purpose, we need to introduce regularity hypotheses on $\partial\Omega$ and some notations. 
We assume that  the boundary~$\partial\Omega$ is smooth enough so that the outward unit normal vector $n(x)$ at $x \in \partial\Omega$ is well-defined as well as $\d\sigma_{\! x}$ the Lebesgue surface measure on $\partial\Omega$. %
The precise regularity on $\partial\Omega$ that we will need is that the signed distance $\delta$ defined by $\delta(x) := - d(x,\partial\Omega)$ if $x \in \Omega$, $\delta(x) :=  d(x,\partial\Omega)$ if $x \in \Omega^c$,
so that $\Omega = \{x \in \R^d, \delta(x) < 0\}$, satisfies $\delta \in W^{3,\infty}(\Omega)$ and $\nabla \delta(x) \ne 0$ for $x \in \partial\Omega$, so that $\nabla \delta / |\nabla \delta|$ coincides with the outward unit normal vector $n$ on $\Omega$. 
We then define $\Sigma_\pm^x := \{ {v} \in \R^d; \pm \, {v} \cdot n(x) > 0 \}$ the sets of outgoing ($\Sigma_+^x$) and incoming ($\Sigma_-^x$) velocities at the point $x \in \partial\Omega$ as well as  
$$
\Sigma_\pm := \Big\{ (x,{v}) \in \Sigma; \pm n(x) \cdot {v} > 0 \Big\} = \Big\{(x,{v}); \, x \in \partial\Omega, \, {v} \in \Sigma^x_\pm \Big \}. 
$$
We denote by $\gamma f$ the trace of $f$ on $\Sigma$, and by $\gamma_{\pm} f = \mathbf 1_{\Sigma_{\pm}} \gamma f$ the traces on $\Sigma_{\pm}$.
The boundary condition \eqref{eq:BdyCond} thus takes into account how particles are reflected by the wall and takes the form of a balance between the values of the trace $\gamma f$  on the outgoing and incoming velocities subsets of the boundary. We assume that the reflection operator acts locally in time and position, namely
$$
(\RRR \gamma_{\!+}  f)(t,x,v) =  \RRR_x (\gamma_{\!+}  f (t,x,\cdot))(v)
$$
and more specifically it is a possibly position dependent Maxwell boundary condition operator 
\beqn\label{eq:boundary} 
\RRR_x (g (x,\cdot))(v)  = (1-\alpha(x))   g (x , R_x v) + \alpha(x) D g  (x,v),
\eeqn
for any $ (x,v) \in \Sigma_-$ and for any function $g : \Sigma_+ \to \R$. Here $\alpha : \partial\Omega \to [0,1]$
is a Lipschitz 
function, called the accommodation coefficient,  $R_x$ is the specular reflection operator
$$
R_x v = v - 2 n(x) (n(x) \cdot v),
$$
and $D$ is the diffusive operator
\begin{align}
\label{eq:def_D}
D g(x,v) = c_\mu \mu(v) \widetilde g (x), \quad \widetilde g (x) = \int_{\Sigma^x_+} g(x,w) \, n(x) \cdot w \, \d w ,
\end{align}
where the constant  $c_\mu := (2\pi)^{1/2}$ is such that  $c_\mu \widetilde{\mu} = 1$ and we recall that $\mu$ stands for the standard Maxwellian \eqref{eq:standardMaxw}. The boundary condition \eqref{eq:boundary} corresponds to the \emph{pure specular reflection} boundary condition when  $\alpha \equiv 0$ 
and it corresponds to the \emph{pure diffusive} boundary condition when $\alpha \equiv 1$. 
It is worth emphasizing that when $\gamma f$ satisfies the boundary condition \eqref{eq:BdyCond}--\eqref{eq:boundary}, for any test function $\varphi = \varphi(v)$ and any $x \in \partial \Omega$, 
\beqn\label{eq:invariantsBoundary}
\int_{\R^d} \gamma f \varphi n(x) \cdot v \, \dv 
=
\int_{\Sigma^x_+} \gamma_{\!+} f  \, n(x) \cdot v \, [ \varphi - (1-\alpha(x)) \varphi \circ R_x - \alpha(x)  c_\mu  \widetilde{\varphi  \circ R_x \mu})] \, \dv. 
\eeqn
As a consequence,  whatever is the accommodation coefficient $\alpha$, making the choice $\varphi = 1$ so that $\varphi \circ R_x = c_\mu \widetilde{\varphi \circ R_x \mu} =1$, we get 
\beqn\label{eq:invariantsBoundary1}
\int_{\R^d} \gamma f \, n(x) \cdot v \, \dv  = 0, 
\eeqn
which means that there is no flux of mass at the boundary (no particle goes out nor enters in the domain). 
Assuming now $\alpha \equiv 0$, making the choice~$\varphi (v) = |v|^2$ and  observing that $|R_x v|^2 = |v|^2$, we get
\beqn\label{eq:invariantsBoundary2}
\int_{\R^d} \gamma f \, |v|^2 \,  n(x) \cdot v \, \dv  = 0, 
\eeqn
which means that  there is no flux of energy at the boundary in the case of the pure specular reflection boundary condition.

\smallskip

\noindent{\it The collisional operator.} Let us now describe the hypotheses made on the collisional linear operator $\CCC$ involved in the linear evolution equation \eqref{eq:dtf=Lf}. 
We assume that the operator acts locally in time and position, namely
$$
(\CCC  f)(t,x,v) = \CCC (  f (t,x,\cdot))(v),
$$
that the operator has mass, velocity and energy conservation laws, namely
\beqn\label{eq:local-conservations}
\int_{\R^d} (\CCC g)(v) \, \varphi (v) \, \dv = 0, 
\eeqn
for $\varphi := 1, v_i, |v|^2$, $i \in \{1,\dots,d\}$, and for any nice enough function $g$, and that the operator has a spectral gap in the classical Hilbert space associated to the standard Maxwellian $\mu$.
In order to be more precise, we introduce the Hilbert space 
$$
L^2_v (\mu^{-1}) := \left\{ f : \R^d \to \R \; \Big| \;  \int_{\R^d} f^2 \mu^{-1} \, \dv < +\infty \right\}  
$$
endowed with the scalar product
$$
(f,g)_{L^2_v(\mu^{-1})} := \int_{\R^d} f g \mu^{-1} \, \dv
$$
and the associated norm $\| \cdot \|_{L^2_v(\mu^{-1})}$. 
We assume that the operator $\CCC$ is a closed operator with dense domain $\mathrm{Dom} (\CCC)$ in $L^2_v (\mu^{-1})$ which satisfies:
\begin{itemize}

\item[(A1)] Its kernel is given by
$$
\mathrm{ker} (\CCC) = \mathrm{span}\{ \mu , v_1 \mu , \ldots , v_d \mu , |v|^2 \mu  \}, 
$$
and we denote by $\pi f$ the projection onto $\mathrm{ker} (\CCC)$ given by 
\beqn\label{def:pi}
\pi f = \left( \int_{\R^d} f \, \d w \right) \mu 
+ \left( \int_{\R^d} w f \, \d w   \right) \cdot v \mu 
+ \left( \int_{\R^d} \frac{|w|^2-d}{\sqrt{2d}} \, f \, \d w   \right) \frac{|v|^2-d}{\sqrt{2d}} \, \mu.
\eeqn

\item[(A2)] The operator is self-adjoint on $L^2_v(\mu^{-1})$ and negative $(\CCC f , f )_{L^2_v(\mu^{-1})} \le 0$, so that its spectrum is included in $\R_{-}$, and \eqref{eq:local-conservations} holds true for any { $g \in \mathrm{Dom} (\CCC)$}. We assume furthermore that $\CCC$ satisfies a coercivity estimate, more precisely that there is a positive constant $\lambda >0$ such that for any $f \in \mathrm{Dom} (\CCC)$ one has
\begin{equation}\label{eq:coercivity}
(-\CCC f , f)_{L^2_v(\mu^{-1})} \ge \lambda \|  f^\perp \|_{L^2_{v}(\mu^{-1})}^2,
\end{equation}
where $f^\perp := f - \pi f$.

\item[(A3)] For any polynomial function $\phi=\phi(v) : \R^d \to \R$ of degree $\le 4$, there holds  $\mu \phi \in \hbox{Dom}(\CCC)$, so that 
there exists a constant $C_\phi \in (0,\infty)$ such that 
$$
\bigl\| \CCC (\phi \mu) \bigr\|_{L^2_v(\mu^{-1})}  \le C_\phi.
$$

\end{itemize}

\subsection{Conservation laws}

Without loss of generality, we shall assume hereafter that the domain $\Omega$ verifies 
\begin{equation}\label{eq:Omega-centre}
|\Omega| = \int_{\Omega} \dx = 1
\quad\text{and}\quad
\int_{\Omega} x \, \dx = 0.
\end{equation}

One easily obtains from \eqref{eq:local-conservations}, the Stokes theorem and   \eqref{eq:invariantsBoundary1} that any solution $f$ to equation \eqref{eq:dtf=Lf}--\eqref{eq:BdyCond} satisfies the conservation of mass
$$
\frac{\d}{\dt} \int_\OO f \, \dv \, \dx 
= \int_\OO (\CCC f - v \cdot \nabla_x f) \, \dv \, \dx
=0.   
$$

In the case of the specular reflection boundary condition, that is \eqref{eq:BdyCond} with $\alpha \equiv 0$, some additional conservation laws appear. On the one hand, one also has the conservation of energy
$$
\frac{\d}{\dt} \int_\OO |v|^2 f \, \dv \, \dx 
= \int_\OO |v|^2(\CCC f - v \cdot \nabla_x f) \, \dv \, \dx
=0,
$$
because of  \eqref{eq:local-conservations}, the Stokes theorem again and   \eqref{eq:invariantsBoundary2}. 
On the other hand, if the domain~$\Omega$ possesses rotational symmetry,   we also have the conservation of the corresponding angular momentum. 
More precisely, we define the set of all infinitesimal rigid displacement fields
\begin{equation}\label{eq:RR}
\RR := \{ x \in \Omega \mapsto Ax + b \in \R^d \,; A  \in \MM^a_{d} (\R), \; b \in \R^d \},
\end{equation}
where $\MM^a_{d} (\R)$ denotes the set of skew-symmetric $d \times d$-matrices with real coefficients,
as well as the linear manifold of \emph{centered} infinitesimal rigid displacement fields preserving $\Omega$
\begin{equation}\label{eq:RROmega}
\RR_\Omega = \{ R \in \RR \mid  b=0 , \;   R(x) \cdot n(x)  = 0, \; \forall \, x \in \partial\Omega  \}.
\end{equation}
We observe here that, thanks to the assumption \eqref{eq:Omega-centre}, we can work only with \emph{centered} infinitesimal rigid displacement fields preserving $\Omega$. Indeed, if $R$ is an infinitesimal rigid displacement field preserving $\Omega$, that is, $R(x) = Ax + b \in \RR$ is such that $R(x) \cdot n(x) = 0$ on $\partial \Omega$, then 
$$
\begin{aligned}
|b|^2 
&= \int_{\Omega} \nabla (b \cdot x) \cdot (Ax + b) \, \dx  \\
&= -\int_{\Omega} (b \cdot x) \Div (Ax + b) \, \dx 
+ \int_{\partial \Omega} (b \cdot x) (Ax+b) \cdot n(x) \, \d\sigma_{\!x} 
= 0,
\end{aligned}
$$
and thus $b=0$. When the set $\RR_\Omega$ is not reduced to $\{ 0 \}$, that is when $\Omega$ has rotational symmetries, then one deduces the conservation of angular momentum 
$$
\frac{\d}{\dt} \int_\OO R(x) \cdot v f \, \dv \, \dx = 0,   \quad \forall \,  R \in \RR_\Omega.
$$
Indeed if $R \in \RR_\Omega$, there exists $A  \in \MM^a_{d} (\R)$ such that $R(x)=Ax$ for any { $x \in \Omega$}. We then compute, using integration by parts,
\begin{align*}
\frac{\d}{\dt} \int_\OO  R(x) \cdot v f \dv \, \dx 
&= \int_\OO Ax \cdot v  (- v \cdot \nabla_{x} f + \CCC f) \, \dv \, \dx \\
&= \int_\OO \partial_{x_k} (Ax \cdot v) v_k  f \, \dv \, \dx 
- \int_\Sigma Ax \cdot v \, { \gamma f } \, n(x) \cdot v \, \dv \, \d\sigma_{\! x} \\ 
&= - \int_\Sigma Ax \cdot v \, { \gamma f } \, n(x) \cdot v \, \dv \, \d\sigma_{\! x}, 
\end{align*}
thanks to the  velocity conservation law \eqref{eq:local-conservations} and the fact that $A$ is skew-symmetric. 
For the boundary term, using \eqref{eq:invariantsBoundary} with $\varphi(x,v) := Ax \cdot v$ and $\alpha \equiv 0$, we get 
\begin{align*} 
\int_\Sigma Ax \cdot v \,\gamma f \, n(x) \cdot v \, \dv \, \d\sigma_{\! x}   
& = \int_{\Sigma_+} Ax \cdot (v - R_x v) \gamma_{+}f \, |n(x) \cdot v| \, \dv \, \d\sigma_{\! x} \\ 
&  = 2\int_{\Sigma_+} (Ax \cdot n(x)) \gamma_{+}f \, |n(x) \cdot v|^2 \, \dv \, \d\sigma_{\! x} =0,
\end{align*}
because $v - R_x v= 2(n(x) \cdot v) n(x)$ and $R \in \RR_{\Omega}$.

\subsection{Main results}\label{ssec:main}

Define the position and velocity dependent Hilbert space
$$
\HH = L^2_{x,v} (\mu^{-1}) 
:= \left\{ f : \OO \to \R \; \Big| \;  \int_\OO f^2 \mu^{-1} \, \dv \, \dx < +\infty \right\} 
$$
endowed with the scalar product
$$
\la f,g \ra_{\HH} := \int_\OO f g \mu^{-1} \, \dv \, \dx
$$
and the associated norm $\| \cdot \|_{\HH}$. 
For $f \in \HH$, we also introduce the following conditions: 
\begin{align}
&\int_\OO f \, \dx \, \dv =0 , \tag{C1}\label{eq:C1}\\ 
&\int_\OO |v|^2 f \, \dx \, \dv = 0 , \tag{C2}\label{eq:C2}\\   
&\int_\OO R(x) \cdot v f \, \dx \, \dv =0, \quad \forall \, R \in \RR_\Omega. \tag{C3}\label{eq:C3}
\end{align}

We are now able to state our main hypocoercivity result:
\begin{theo}\label{theo:hypo}
There exists a scalar product $\la\!\la \cdot , \cdot \ra \! \ra$ on the space $\HH$ so that the associated norm $\Nt \cdot \Nt$ is equivalent to the usual norm $\| \cdot \|_{\HH}$, and for which the linear operator $\LLL$ satisfies the following coercivity estimate: there is a positive constant $\kappa >0$ such that 
$$
\la \! \la - \LLL f , f \ra\!\ra \ge \kappa  \Nt f \Nt^2
$$
for any $f \in \mathrm{Dom}(\LLL)$ satisfying the boundary condition \eqref{eq:BdyCond}, assumption~\eqref{eq:C1} and furthermore assumptions \eqref{eq:C2}-\eqref{eq:C3} in the specular reflection case ($\alpha\equiv 0$ in \eqref{eq:BdyCond}).
\end{theo}

\noindent This result improves existing results regarding hypocoercivity in a bounded domain for the {linearized} Boltzmann and Landau equations (and consequently for their long-time stability, see Theorem~\ref{theo:main}) in three regards:

\smallskip
-- We consider a general, smooth enough, convex or non-convex domain.

\smallskip
-- The $L^2$ estimates that we establish are constructive, which means that they depend constructively of some collisional constants (that appear in the estimates (A2)-(A3) satisfied by the collisional operator $\CCC$) and some geometrical constants depending on the domain $\Omega$ (that appear in some Poincaré and Korn inequalities which can be made explicit, at least for a domain with simple geometry). 

\smallskip
-- Our method encompasses the three boundary conditions (pure diffusive, specular reflection and Maxwell) in a single treatment. In particular, we can solve the Maxwell boundary condition in the case where the accommodation coefficient $\alpha$ vanishes everywhere or on some subset of the boundary. 

\medskip

Our proof is based on a $L^2$-hypocoercivity approach. The challenge of hypocoercivity is to understand the interplay between the collision operator that provides dissipativity in the velocity variable and the transport one which is conservative, in order to obtain global dissipativity for the whole problem. There are two main hypocoercivity methods, the $H^1$ and the $L^2$ ones. The $H^1$-hypocoercivity approach has been first introduced for hypoelliptic operators by H\'erau, Nier~\cite{Herau-Nier} and Eckmann, Hairer~\cite{MR1969727},  further developed by Nier, Helffer~\cite{HN05} and Villani~\cite{Mem-villani} and extended to more general kinetic operators in Villani~\cite{Mem-villani} and Mouhot, Neumann~\cite{Mouhot-Neumann}. It is also reminiscent of the work by Desvillettes and Villani on the trend to global equilibrium for spatially inhomogeneous kinetic systems in \cite{Desvillettes-Villani-2001}, \cite{Desvillettes-Villani-2005}, and of the high order Sobolev energy method developed by Guo in \cite{Guo-2002-I} and subsequently. 
In summary, the idea consists in endowing the~$H^1$~space with a new scalar product which makes coercive  the considered operator and whose associated norm is equivalent to the usual $H^1$ norm.  
In order to be adapted to more general operators and geometries, the $L^2$-hypocoercivity technique for one dimensional space of collisional invariants has been next introduced by H\'erau~\cite{MR2215889}  and developed by Dolbeault-Mouhot-Schmeiser~\cite{Dolbeault2009511,MR3324910}. The $L^2$-hypocoercivity technique for a space of collisional invariants of dimension larger than one (including the Boltzmann and Landau cases) has been introduced by Guo in \cite{MR2679358}, and developed further mainly by Guo, collaborators and students. Again the idea consists in endowing the $L^2$ space with a new scalar product which makes coercive  the considered operator and whose associated norm is equivalent to the usual $L^2$ norm. 
 
\medskip
 
We present hereafter the line of reasoning of this last approach that will be ours. It heavily relies on the micro-macro decomposition of the solution of the equation: $f = f^\perp + \pi f$, where~$f^\perp$ denotes the microscopic part and $\pi f$ the macroscopic part defined in~\eqref{def:pi}. The coercive estimate~\eqref{eq:coercivity} on the collision operator~$\CCC$ already gives a control on~$f^\perp$ but not on the macroscopic term $\pi f$. Then, in order to control the macroscopic part, we construct a new scalar product on~$\HH$ by adding, step by step, new terms in order to control the missing terms appearing on the macroscopic part~$\pi f$. Roughly speaking, the scalar product that we cook up takes the following form:
\bean
\la \! \la f, g \ra \! \ra
&:=& \la f,g \ra_{\HH} - \eta \la \widetilde\pi f , \nabla \Delta^{-1}\pi g \ra_{L^2_x(\Omega)}  - \eta \la \nabla \Delta^{-1} \pi  f ,  \widetilde\pi g \ra_{L^2_x(\Omega)} , 
\eean
choosing $\eta > 0$ small enough, and where the moments operator { $ \widetilde\pi : \HH \to (L^2_x(\Omega))^{d}$} and the inverse Laplacian type operator $\Delta^{-1}$ have to be suitably defined (see Sections~\ref{sec:elliptic}~\&~\ref{sec:proof}).

\smallskip
Our proof is a variant of previous proofs of the same type but differs from them by several aspects: 

\smallskip
(i) The order between the $\nabla$ operator and the $\Delta^{-1}$ operator is the one from Guo's approach~\cite{MR2679358,MR3562318} rather than the one from Dolbeault-Mouhot-Schmeiser's approach~\cite{Dolbeault2009511,MR3324910}. That is important in order to handle the rather singular operator involved by the boundary condition. 

\smallskip
(ii) The choice of the mean operator $\widetilde \pi f$ differs from the one used in \cite{MR2679358,MR3562318,MR3579575,MR3840911} but looks very much like the one in \cite{MR2813582,MR2966364,CDHMMS}. It allows to deal with 
general Maxwell boundary condition (and the possibility that $\alpha$ vanishes somewhere or everywhere) but leads to a first natural control of the symmetric gradient of the momentum component of the macroscopic part~$\nabla^s m$ instead of the full derivative $\nabla m$ as in Guo's approach.
 
\smallskip
(iii) The definition of the $\Delta^{-1}$ operator has to be chosen wisely in order to handle the general Maxwell boundary condition and the mean operator $\widetilde \pi f$. 
We thus need to establish natural $H^{-1} \to H^1$ and $L^2 \to H^2$ regularity estimates for some classical elliptic problems but associated with somehow unusual boundary conditions. 

\smallskip
Let us give a few more details about (iii). First, we shall introduce an auxiliary Poisson equation with Robin or Neumann boundary conditions, which are devised in order to control \emph{mass} and \emph{energy} terms of $\pi f$. This result is stated in Theorem~\ref{theo:Poisson} and is based on Poincar\'e type inequalities.
Next, we shall introduce a tailored Lam\'e-type system with mixed Robin-type boundary conditions in order to deal with the \emph{momentum} component of the macroscopic part $\pi f$. The corresponding result is presented in Theorem~\ref{theo:regH2-korn} and is based on Korn-type inequalities, which are discussed in Section~\ref{sec:Korn}. For more information on Korn inequalities we refer to the fundamental result of Duvaut-Lions \cite[Theorem~3.2 Chap.~3]{MR0521262}, and on the variant introduced by Desvillettes and Villani \cite{DV02}. For further references and a recent treatment of Korn's inequality, we refer to Ciarlet and Ciarlet \cite{MR2119999}.
For more details concerning the regularity 
issue for similar elliptic equations and systems we refer to~\cite{MR775683,CostabelDN,MR1452171} and the references therein.

\medskip
Let us now point out that our hypocoercivity result obtained in Theorem~\ref{theo:hypo} enables us to deduce an exponential stability result for our equation~\eqref{eq:dtf=Lf} supplemented with the boundary condition~\eqref{eq:BdyCond}.

\begin{theo}\label{theo:main}
Let $f_{\mathrm{in}} \in \HH$ satisfying assumption~\eqref{eq:C1} and furthermore assumptions~\eqref{eq:C2} and \eqref{eq:C3} in the specular reflection case ($\alpha\equiv 0$ in \eqref{eq:BdyCond}).
There exist positive constants $\kappa , C >0$ such that for any solution $f$ to \eqref{eq:dtf=Lf}--\eqref{eq:BdyCond} associated to the initial data~$f_{\mathrm{in}}$, there holds
$$
\| f(t) \|_{\HH} \le C e^{-\kappa t} \| f_{\mathrm{in}} \|_{\HH}, \quad \forall \, t \ge 0.
$$
\end{theo}

\noindent This result is a first step towards  the global existence and the study of the long-time behavior of solutions to the nonlinear problem~\eqref{eq:nonlinear}-\eqref{eq:Bdynonlinear} in a close-to-equilibrium regime that will be the object of a forthcoming work. 

\smallskip
We here briefly mention some similar coercivity estimates or exponential stability results established in the last decade for linear kinetic equations (mainly for the linearized Boltzmann equation) in a bounded domain.
These ones have then been used for proving global existence of solutions to nonlinear equation in a close-to-equilibrium regime and convergence to the equilibrium in the long-time asymptotic. 
As already mentioned, Guo~\cite{MR2679358} has first proved a $L^2_{x,v}$ coercivity estimate for the cutoff Boltzmann equation with hard potentials or hard-spheres by using non-constructive technique in two cases:  
the specular reflection boundary condition with strictly convex and analytic domains~$\Omega$ and the pure diffusive boundary condition assuming the domain $\Omega$ is smooth and convex. 
These results have been generalized by  Briant and Guo \cite{MR3562318} who derived constructive  exponential stability estimates in $L^2_{x,v}$ for any positive and constant accommodation coefficient $\alpha \in (0,1)$,  with no more convexity assumptions on $\Omega$. For the same equation endowed with specular reflection boundary condition, a still non-constructive $L^2$ estimate was derived in the convex setting, without analyticity assumptions on the domain, by Kim and Lee \cite{MR3762275}.  The authors then extended their results to periodic cylindrical domain with non-convex analytic cross-section~\cite{MR3840911}.

Furthermore, the only results we are aware of in the case of long-range interaction, that is, for non-cutoff Boltzmann and Landau collision operators in a bounded domain, are the very recent works of Guo-Hwang-Jang-Ouyang~\cite{MR4076068} (see also \cite{Guoetal-specular2}) for the Landau equation with specular reflection boundary condition, and Duan-Liu-Sakamoto-Strain~\cite{DuanLiuSakamotoStrain} for non-cutoff Boltzmann and Landau equations in a finite channel with inflow or specular reflection boundary conditions. However, as far as we understand, the arguments presented in 
\cite{MR4076068} seem to be constructive only when $\partial\Omega$ is flat, while the arguments presented in~\cite{Guoetal-specular2} are again non-constructive. 

It is also worth mentioning that an alternative existence of solutions framework to the above quite strong but close-to-equilibrium regime framework has been introduced by DiPerna and Lions who proved  in \cite{DiPernaLionsAnnals,MR1088276,MR1295942} the existence of global weak (renormalized) solutions of arbitrary amplitude to the Boltzmann equation  in the case of the whole space 
for initial data satisfying only the physically natural condition that the total mass, energy and entropy are finite.
The extension to the case of  a bounded domain with reflection conditions (including specular reflection, pure diffusive reflection and Maxwell reflection) has been then obtained  in  \cite{HamdacheARMA1992,ArkerydMaslova,MischlerCMP2000,MischlerENS2010}. 
We must emphasize that our treatment of boundary terms bears some similarity with the analysis made in \cite{MischlerENS2010} in order to   take advantage of the information provided by Darroz\`es and Guiraud inequality~\cite{DGineq}.

\medskip

To end this introduction, we point out that in Section~\ref{sec:weak}, we broaden our study to the case where the linearized operator only enjoy a weak coercivity estimate to obtain results of weak hypocoercivity and sub-exponential stability in Theorems~\ref{theo:weak-hypo} and~\ref{theo:weak-main}.

\smallskip

Also, in Section~\ref{sec:hydro}, we extend our study to a rescaled version of~\eqref{eq:dtf=Lf} which naturally arises in the analysis of hydrodynamical limit problems, we obtain hypocoercivity and stability results uniformly with respect to the rescaling parameter in Theorems~\ref{theo:hydro-hypo} and~\ref{theo:hydro-main}.

\medskip  
\noindent {\bf Acknowledgements.} The authors thank O.~Kavian and F.~Murat for enlightening  discussions and for having pointing out several relevant references.
This work has been partially supported by the Projects EFI: ANR-17-CE40-0030 (K.C.\ and I.T.) and SALVE: ANR-19-CE40-0004 (I.T.) of the French National Research Agency (ANR). A.B. acknowledges financial support from R\'egion \^Ile de France.


\section{Elliptic equations}
\label{sec:elliptic}

We present some functional estimates associated to some elliptic problems related to the macroscopic quantities. 
In this section, we denote the classical norm on $L^2_x(\Omega)$ by $\|\cdot\|$ and  the associated scalar product by~$(\cdot,\cdot)$. We also write 
$$
\langle f \rangle :=  \int_\Omega f \, \dx
$$ the mean of $f$ (recall our normalization assumption \eqref{eq:Omega-centre}). The operators that we consider only act on the position variable $x$, so that, in order to lighten the notations, we will not mention it in our proofs. For the same reason, we often write $\partial_i$ for $\partial_{x_i}$, $i \in \{1,\dots, d\}$.

\subsection{Poincar\'e inequalities and Poisson equation}
\label{sec:poincare}

We consider the following Poisson equation
\begin{equation}\label{eq:elliptic}
\left\{
\begin{aligned}
- \Delta u &= \xi \quad\text{in}\quad \Omega , \\
(2-\alpha(x))\nabla u \cdot n(x) + \alpha(x) u &= 0 \quad\text{on}\quad \partial \Omega,   
\end{aligned}
\right.
\end{equation}
for a scalar source term $\xi : \Omega \to \R$. 
Remark that when $\alpha \equiv 0$ then \eqref{eq:elliptic} corresponds to the Poisson equation with homogeneous Neumann boundary condition. Otherwise,~\eqref{eq:elliptic} corresponds to the Poisson equation with homogeneous Robin (or mixed) boundary condition.

We define the Hilbert spaces 
$$
V_1 :=  H^1(\Omega) \quad \hbox{and} \quad
V_0 := \left\{ u \in H^1(\Omega); \ \int_\Omega u \, \dx = 0   \right\}
$$
endowed with the $H^1(\Omega)$-norm, and next
$$
V_\alpha := 
\begin{cases}
V_1 \quad\text{if}\quad \alpha \not\equiv 0 \\
V_0 \quad\text{if}\quad \alpha \equiv 0.
\end{cases}
$$
On $V_\alpha$, we define  the bilinear form 
$$
a_\alpha(u,v) :=   \int_\Omega \nabla u \cdot \nabla v \, \dx 
+ \int_{\partial\Omega}  \frac{\alpha }{2-\alpha} \, u v \, \d\sigma_{\! x}.
$$
We start by a result on Poincaré-type inequalities:
\begin{prop}
There hold
\beqn\label{eq:PWineq}
\forall \, u \in V_0, \quad \| u \| \lesssim \| \nabla u \|, 
\eeqn
and 
\beqn\label{eq:PTypeineq}
 \forall \, u \in V_1, \quad \| u  \|^2 \lesssim a_\alpha(u,u).
\eeqn
\end{prop}

The first inequality is nothing but the classical Poincaré-Wirtinger inequality. For the second inequality (which is probably also classical), we have no precise reference for a constructive proof. For the sake of completeness and because we will need to repeat that kind of argument in the next section, we give a sketch of a non constructive proof by contradiction based on a compactness argument. 

\medskip
\noindent {\it Proof of~\eqref{eq:PTypeineq}.} 
Assuming that \eqref{eq:PTypeineq} is not true, there exists a sequence $(u_n)_{n \in \N}$ in $H^1(\Omega)$ such that 
$$
1 = \| u_n   \|^2 \ge n \left( \| \nabla u_n \|^2 + \left\|  \sqrt{\frac{\alpha}{2-\alpha}} u_n  \right\|_{L^2(\partial\Omega)}^2 \right). 
$$
As a consequence, up to the extraction of a subsequence, there exists  $u \in H^1(\Omega)$ such that $u_n \wto u$ weakly in $H^1(\Omega)$ and $u_n \to u$ strongly in $L^2(\Omega)$. From the above estimate we deduce that $\|  \nabla u  \| \le \liminf_{n \to \infty} \| \nabla u_n  \| = 0$, so that $u = C$ is a constant. On the one hand, we have $\| \sqrt{\alpha/(2-\alpha)} u  \|_{L^2(\partial\Omega)} = \lim_{n \to \infty} \| \sqrt{\alpha/(2-\alpha)} u_n  \|_{L^2(\partial\Omega)} = 0$ so that~$C = 0$. On the other hand, we get $\|  u  \| = \lim_{n\to \infty} \| u_n  \| = 1$, which implies that $C \not = 0$ and thus a contradiction. 
\qed

\smallskip

We now state a result on the existence, uniqueness and regularity of solutions to \eqref{eq:elliptic}.

\begin{theo}\label{theo:Poisson}
For any given $\xi \in L^2(\Omega)$, 
there exists a unique $u \in V_\alpha$ solution to the variational problem
\beqn\label{eq:varValpha}
a_\alpha(u,w)  = ( \xi , w ) ,\quad\forall \, w \in V_\alpha. 
\eeqn
Assuming furthermore that $\la \xi \ra = 0$ when $\alpha \equiv 0$, there holds $u \in H^2(\Omega)$,  $u$ verifies the elliptic equation \eqref{eq:elliptic} a.e.\ and 
\beqn\label{eq:varValphaH2}
\| u \|_{H^2(\Omega)} \lesssim \| \xi \|.
\eeqn
\end{theo}

We give a sketch of the proof of Theorem~\ref{theo:Poisson} which is very classical, except maybe the way we handle the $H^2$ regularity estimate. The proof will be taken up again in the next section where we deal with an elliptic system of equations associated to the symmetric gradient. 

\begin{proof}[Proof of Theorem~\ref{theo:Poisson}] 
We split the proof into 4 steps. The first one is dedicated to the application of Lax-Milgram theorem. The last three ones are devoted to the proof of the $H^2$ regularity estimate: in Step 2, we develop a formal argument which leads to a directional regularity estimate supposing that the variational solution~$u$ is {\em a~priori} smooth; we then make it rigorous in Step 3 by not supposing any smoothness assumption on $u$ and in Step~4, we end the proof of~\eqref{eq:varValphaH2}.

\smallskip\noindent{\sl Step 1. }
We first observe that there exists $\lambda > 0$ such that 
$$
a_\alpha(u,u) \ge \lambda \| u \|_{H^1(\Omega)}^2, \quad \forall \, u \in V_\alpha,
$$
and thus $a_\alpha$ is coercive. The above estimate is a direct consequence of the Poincar\'e-Wirtinger inequality~\eqref{eq:PWineq} in the case when $\alpha \equiv 0$ and the variant of the classical Poincar\'e inequality given in~\eqref{eq:PTypeineq} when $\alpha \not\equiv 0$.  Because $\xi \in L^2(\Omega) \subset V_\alpha'$, we may use the Lax-Milgram theorem and we get the existence and uniqueness of $u \in V_\alpha$ satisfying \eqref{eq:varValpha} as well as 
\beqn \label{eq:u_H1}
\|u\|_{H^1(\Omega)} \lesssim \|\xi\|.
\eeqn

For the remainder of the proof, we furthermore assume  $\la \xi \ra = 0$ when $\alpha \equiv 0$. We claim that \eqref{eq:varValpha} can be improved into the following new formulation: there exists a unique $u \in V_\alpha$ satisfying
\begin{equation}\label{eq:varValpha-bis}
a_\alpha(u,w) = ( \xi , w ) , \quad \forall \, w \in H^1(\Omega).
\end{equation}
When $\alpha \not\equiv 0$ formulation \eqref{eq:varValpha-bis} is nothing but \eqref{eq:varValpha}. 
In the case $\alpha \equiv 0$ so that $V_\alpha \not= H^1(\Omega)$, we remark that for any $w \in H^1(\Omega)$, we have $w - \la w \ra \in V_0$ and therefore
$$
\begin{aligned}
a_\alpha (u,w) 
&= a_\alpha (u , w - \la w \ra) \\
&= \int_{\Omega} \xi  w \, \dx  
- \int_{\Omega} \xi \la w \ra   \dx  = \int_{\Omega} \xi w \, \dx ,
\end{aligned}
$$
where we have used the formulation \eqref{eq:varValpha} and the condition $\la \xi \ra =0$ so that $\int_{\Omega} \xi \la w \ra   \dx = 0$  in the second line.

\smallskip\noindent{\sl Step 2. A priori directional estimate.}
For any small enough open set $\omega \subset \Omega$, we fix 
a vector field $a \in C^2(\bar\Omega)$ such that $|a| = 1$ on $\omega$ and $a \cdot n = 0$ on $\partial \Omega$, 
and we set $X := a \cdot \nabla$ the associated differential operator.
For a smooth function $u$, we compute
\bean
\| \nabla X u\|^2 
 &=&  (\nabla  u, X^* \nabla  X u  ) +  ([\nabla,X]  u,   \nabla  X u  )  
\\
&=&
 (\nabla u,   \nabla  X^*X u  ) +  (\nabla u, [X^*, \nabla ] X u  ) + ([\nabla,X]  u,   \nabla  X u  ),
\eean
where we have used that 
\beqn \label{eq:X*}
(Xf,g) = (f,X^* g), \quad X^* g :=  -  \hbox{div} (a g), 
\eeqn
because  $a \cdot n = 0$ on $\partial \Omega$. On the other hand, we compute formally
\beqn \label{eq:IPPbord}
\int_{\partial\Omega} (Xu)^2 \, \frac{\alpha}{2-\alpha}  \d\sigma_{\! x}
=  \int_{\partial\Omega}  \frac{\alpha}{2-\alpha} \, u ( X^* Xu) \, \d\sigma_{\! x} 
- \int_{\partial\Omega}   \biggl( X \frac{\alpha}{2-\alpha} \biggr) \, u (Xu) \, \d\sigma_{\! x} .
\eeqn
In the next step of the proof, we will work with a discrete version of the operator~$X$ which will allow us to make rigorous computations.
Assuming furthermore now that $u \in V_\alpha$ satisfies \eqref{eq:varValpha-bis} and that $X^*X u \in H^1(\Omega)$, we may use \eqref{eq:varValpha-bis} with { $w := X^*X u $} and we deduce  
\bean
&&\| \nabla X u\|^2 + \int_{\partial\Omega} \frac{\alpha}{2-\alpha} (Xu)^2 \, \d\sigma_{\! x}  
\\
&&\quad=
 (\xi, X^*X u  )  +  (\nabla u, [X^*, \nabla ] X u  ) + ([\nabla,X]  u,   \nabla  X u  )  
  - \int_{\partial\Omega}   \biggl( X \frac{\alpha}{2-\alpha} \biggr) \, u (Xu) \, \d\sigma_{\! x}. 
\eean
We easily compute for $i=1,\dots,d$
\bean
[\partial_i,X] = (\partial_i a) \cdot \nabla,  \quad {[{X^*}, {\partial_i}]} =  \partial_i (\hbox{div} a) + (\partial_i a) \cdot \nabla,
\eean
so that for some constant $C = C(\| a \|_{W^{2,\infty}(\Omega)})$ and any function { $w \in H^1(\Omega)$ }, we have 
$$
\| [\nabla,X]  w \| \le C \| \nabla w \|, \quad
\|  [X^*, \nabla ]  w \| \le C \|  w \|_{H^1(\Omega)}.
$$
We then deduce that for some constant $C = C( \| a \|_{W^{2,\infty}(\Omega)},\|\alpha\|_{W^{1,\infty}(\Omega)})$, we have
\bean
&\| \nabla X u\|^2 
\le
\| \xi \| \|  X^*X u  \| + C \| \nabla u \| \| X u \|_{H^1(\Omega)}  \\
&\hspace{5cm}  +C \| \nabla u \| \|  \nabla  X u \| + C \| u \|_{L^2(\partial\Omega)} \| X u \|_{L^2(\partial\Omega)}.
\eean
Recalling \eqref{eq:u_H1} and observing that $ \|  X^* w \| +  \|  X w \| + \| w \|_{L^2(\partial\Omega)} \lesssim \| w \|_{H^1(\Omega)}$, we obtain 
\bean
\| \nabla X u\|^2 
\lesssim
\| \xi \| \| \nabla X u\| +   \|\xi\|^2,
\eean
and we conclude that
\beqn\label{eq:RegDirXu}
\| \nabla X u\|   \lesssim  \| \xi \| .
\eeqn

\smallskip\noindent{\sl Step 3. Rigorous directional estimate.}
When we do not deal with an {\em a priori} smooth solution, but just with a variational solution $ u \in V_\alpha$ satisfying \eqref{eq:varValpha-bis},
we have to modify the argument in the following way. We define $\Phi_t : \bar\Omega \to \bar\Omega$ the flow associated to 
the differential equation 
\beqn \label{eq:flow}
\dot y = a(y), \quad y(0) = x, 
\eeqn
so that $\Phi_t(x) := y(t)$, $(t,x) \mapsto \Phi_t(x)$ is $C^1$ and $\Phi_t$ is a diffeomorphism on both $\Omega$ and~$\partial\Omega$ for any $t \in \R$.
  We next define 
$$
X^hu(x) := \frac{1}{h} \bigl( u(\Phi_h(x)) - u(x)),
$$
so that $X^h u \in H^1(\Omega)$ if $u \in  V_\alpha$.  Repeating the argument of Step 1, we get the identity 
\beqn \label{eq:nablaXh}
\begin{aligned}
&\| \nabla X^h u\|^2 + \int_{\partial\Omega} \frac{\alpha}{2-\alpha} (X^h u)^2 \, \d\sigma_{\! x} = (\xi, X^{h*}X^h u  )  
+   (\nabla u, [X^{h*}, \nabla ] X^h u  )  
\\
&\qquad \qquad \qquad  + \, \, ([\nabla,X^h]  u,   \nabla  X^h u  )  - \int_{\partial\Omega} u (\Phi_h(x)) \biggl((X^h u) \, X^h \bigg(\frac{\alpha}{2-\alpha}\bigg) \biggr) (x)   \, \d\sigma_{\! x}  ,
\end{aligned}
\eeqn
where we denote
\bean
X^{h*} w (x) := \frac{1}{h} \Big[ w(\Phi_{-h}(x)) \left|\hbox{\rm det} D\Phi_{-h}(x)\right|- w(x) \Big].
\eean
Notice here that we used a discrete version of the integration by parts leading to~\eqref{eq:IPPbord} and it only relies on a change of variable on $\partial \Omega$, which makes our computation fully rigorous.
As in the second step of the proof, we are now going to bound each term of the right-hand-side of~\eqref{eq:nablaXh}. First, notice that for $|h| \le 1$, we have for some $|h_0| \le 1$:
\[
X^h u(x) =  \sum_j \partial_j u(\Phi_{h_0}(x)) a_j(\Phi_{h_0}(x)) 
\]
so that there exists $C=C(\|a\|_{W^{1,\infty}(\Omega)})$ such that for any $|h| \le 1$, we have
$\|X^h u\| \le C \|\nabla u\|$.
We can estimate $\|X^{h*}w\|$ in a similar way using that
$$
X^{h*} w (x) = {\frac1{h}} \big[ w(\Phi_{-h}(x)) - w(x) \big] \left|\hbox{\rm det} D\Phi_{-h}(x)\right|+ {\frac1{h}}  w(x) \big[\left|\hbox{\rm det} D\Phi_{-h}(x)\right| -\left|\hbox{\rm det} D\Phi_{0}(x)\right| \big].
$$
Consequently, we deduce that there exists $C = C(\|a\|_{W^{2,\infty}(\Omega)})$ such that for $|h| \le 1$, 
\beqn
\label{eq:control_Xh}
\|X^{h*} w\| + \|X^h w\| + \|w\|_{L^2(\partial \Omega)} \le C \|w \|_{H^1(\Omega)}.
\eeqn
For $i=1,\dots,d$, { for $|h| \le 1$, $x \in \bar{\Omega}$, writing $\Phi_h(x) = (\Phi_{h,1}(x),\dots,\Phi_{h,d}(x))$}, we compute
\begin{align*}
[\partial_i, X^h] w (x)
&= {\frac1{h}} \sum_{j \neq i} \partial_j w (\Phi_h(x)) \partial_i \Phi_{h,j}(x) + {\frac1{h}} \partial_i w (\Phi_h(x)) \left( \partial_i \Phi_{h,i}(x) - 1\right) \\
&= {\frac1{h}} \sum_{j} \partial_j w (\Phi_h(x)) \left(\partial_i \Phi_{h,j}(x) - \partial_i \Phi_{0,j}(x)\right)
\end{align*}
and similarly
\begin{align*}
[X^{h*}, \partial_i ] w (x)
&= {1 \over h} \sum_{j} \partial_j w (\Phi_{-h}(x)) \left(\partial_i \Phi_{0,j}(x)-\partial_i \Phi_{-h,j}(x) \right) \left|\hbox{\rm det} D\Phi_{-h}(x)\right| \\
&\quad - {1 \over h} w(\Phi_{-h}(x)) \partial_i \left|\hbox{\rm det} D\Phi_{-h}(x)\right|.
\end{align*}
As previously, we can easily bound $[\partial_i, X^h] w$ and the first term in $[X^{h*}, \partial_i ] w$ by $C \|\nabla w\|$ with $C=C(\|a\|_{W^{1,\infty}(\Omega)})$ for any $|h| \le 1$. The second term of $[X^{h*}, \partial_i ] w$ can be bounded by $C \|w\|$ with $C=C(\|a\|_{W^{2,\infty}(\Omega)})$ for any $|h| \le 1$ since for any $j$, we have $\partial_{ij} \Phi_0(x) = 0$.
This implies that there exists $C = C(\| a \|_{W^{2,\infty}(\Omega)})$ such that for $|h| \le 1$ and any function $w$ in $H^1(\Omega)$, we have 
$$
\| [\partial_i,X^h]  w \| \le C \| \nabla w \|, \quad
\|  [X^{h*}, \partial_i ]  w \| \le C \|  w \|_{H^1(\Omega)}.
$$
We deduce that for some $C = C(\|a\|_{W^{2,\infty}}, \|\alpha\|_{W^{1,\infty}})$, we have for any $|h| \le 1$:
\begin{align*}
\| \nabla X^h u\|^2 
&\le
\| \xi \| \|  X^{h*} X^h u  \| + C \| \nabla u \| \| X^h u \|_{H^1(\Omega)}  \\
&\qquad +C \| \nabla u \| \|  \nabla  X^h u \| + C \| u \|_{L^2(\partial\Omega)} \| X^h u \|_{L^2(\partial\Omega)}
\end{align*}
and then, { using $\|X^h u \| \lesssim \| \nabla u\|$, \eqref{eq:control_Xh} and \eqref{eq:u_H1},} 
$$
\| \nabla X^h u\|   \lesssim  \| \xi \| .
$$
Passing to the limit $h\to0$, we recover \eqref{eq:RegDirXu}. 
 
\smallskip\noindent
{\sl Step 4. Proof of~\eqref{eq:varValphaH2}.} 
Consider a small enough open set $\omega \subset \Omega$, so that we may fix $a^1, \dots , a^d$ a family of smooth vector fields such that it is an orthonormal basis of $\R^d$ at any point $x \in \omega$ and~$a^1(x) = n(x)$ for any $x \in \partial\Omega \cap \partial\omega$. In order to see that it indeed holds true, we may argue as follows. If $ \partial\Omega \cap \partial\omega = \emptyset$, we may take $a^j := e_j$ the canonical basis of $\R^d$. Otherwise, we fix $x_0 \in  \partial\Omega \cap \partial\omega$. Because 
$\nabla \delta(x_0) \ne 0$, 
we may fix first $i \in \{1, \dots, d \}$ such that $\partial_{x_i} \delta(x_0) \ne 0$ and thus $\partial_{x_i} \delta(x)  \ne 0$
for any $x \in \omega$, for $\omega$ small enough. We then define $b^1:= \nabla \delta$, $b^j := e_{j-1}$ for any $j \in \{2, \dots, i \}$ and $b^j := e_{j}$ for any $j \in \{i+1, \dots, d \}$.
Finally, we apply the Gram-Schmidt process to $(b^1(x),\dots,b^d(x))$ to obtain $(a^1(x),\dots,a^d(x))$. 
We set now $X_i := a^i \cdot \nabla$. From the third step, we have 
\beqn\label{eq:RegDirXiu}
\| \nabla X_i u\|   \lesssim  \| \xi \|, \quad \forall \, i = 2, \dots, d.
\eeqn
As a consequence of our previous construction, the matrix $A:=(a^1,\dots, a^d)$ is orthonormal. 
We thus have $\delta_{k\ell} = a^k \cdot a^\ell = a_k \cdot a_\ell$, where we denoted by $a_m$ the $m$-th line vector of the matrix~$A$. 
As a consequence, we have
\beqn \label{eq:Deltau}
\sum_i X^*_i X_i u 
= - \sum_{i,k,\ell} \partial_k (a^i_k a^i_\ell \partial_\ell u) = - \sum_{k,\ell} \partial_k(a_k \cdot a_\ell \partial_\ell u) 
=
- \Delta u, 
\eeqn
from what we deduce
\bean
{ X^*_1 X_1 u 
= \xi -  \sum_{i \neq 1} X^*_i X_i u. 
}
\eean
Because of \eqref{eq:RegDirXiu}, the above identity and $[X_1^*,X_1] u = {(a^1 \cdot \nabla \hbox{div} (a^1))u }$, we get 
\begin{align*}
\|X_1^2 u\|^2 &= (X_1^* X_1 u, X_1^* X_1 u) + (X_1 u, [X_1^*,X_1] X_1 u) \\
&\lesssim \|\xi\|^2 + \sum_{i \ne 1} \Big( \| \nabla X_i u\| \|\xi\| + \| \nabla X_i u\|^2 \Big) + \|u\|_{H^1(\Omega)}^2 \lesssim \|\xi\|^2. 
\end{align*} 
Together with \eqref{eq:RegDirXiu} again, we have then established
 \beqn\label{eq:RegDirXiju}
\| X_i X_j u\|   \lesssim  \| \xi \|, \quad \forall \, i ,j = 1, \dots, d.
\eeqn
Recalling that $A = (a^1, \dots, a^d)$, we have $\partial_i = (A X)_i$.   As a consequence, we may write 
\bean
\partial_i\partial_j u 
&=&  \sum_{m,\ell}  A_{im} X_m A_{j\ell} X_\ell u 
\\
&=&  \sum_{m,\ell} \left(A_{im}  A_{j\ell} X_m X_\ell u  +  A_{im}  [ X_m , A_{j\ell} ] X_\ell u\right),  
\eean
where the last operator is of order $1$. Together with the starting point estimate~\eqref{eq:u_H1} and~\eqref{eq:RegDirXiju}, we conclude that 
$$
\| \partial_i \partial_j u\|   \lesssim  \| \xi \|, \quad \forall \, i ,j = 1, \dots, d, 
$$
which ends the proof of~\eqref{eq:varValphaH2}. We can now conclude the proof of Theorem~\ref{theo:Poisson}. Indeed, because~$u \in H^2(\Omega)$, we may compute from~\eqref{eq:varValpha} and the Stokes formula:
\bean
\int_{\partial\Omega} \left\{ \frac{\partial u}{\partial n} + \frac{\alpha u}{2-\alpha}  \right\} w \, \d\sigma_{\! x} 
&=& \int_\Omega  \Delta u \, w \, \dx +  \int_\Omega  \nabla u \cdot \nabla w \, \dx  + \int_{\partial\Omega} \frac{\alpha }{2-\alpha} \, u  w \, \d\sigma_{\! x} 
\\
&=& \int_\Omega  (\Delta u  + \xi ) w \, \dx,
\eean
for any $w \in V_\alpha$. Considering first $w \in C^1_c(\Omega)$ and next $w \in C^1(\bar\Omega)$, we get that $u$ satisfies both equations in  \eqref{eq:elliptic}.
\end{proof}

\subsection{Korn inequalities and the associated elliptic equation}
\label{sec:Korn}

For a vector field $M=(m_i)_{1 \le i \le d} : \Omega \to \R^d$, we define its symmetric gradient through
$$
\nabla^{s}_x M := {1 \over 2} \left( \partial_j{m_i} + \partial_i m_j \right)_{1 \le i,j \le d},  
$$
as well as its skew-symmetric gradient by 
$$
\nabla^{a}_x M := {1 \over 2} \left( \partial_j{m_i} - \partial_i m_j \right)_{1 \le i,j \le d}.  
$$
Through this section, in order to lighten the notations, we will write $\nabla^s$ for $\nabla^s_x$, and~$\nabla^a$ for $\nabla^a_x$.
We consider the system of equations 
\begin{equation}\label{eq:elliptic-korn}
\left\{
\begin{aligned}
- \Div (\nabla^s U)  = \Xi \quad&\text{in}\quad \Omega , \\
U \cdot n = 0 \quad&\text{on}\quad \partial \Omega, \\  
(2-\alpha) \left[\nabla^s U  n - (\nabla^s U : n \otimes n) n \right] + \alpha U = 0 \quad&\text{on}\quad \partial \Omega,
\end{aligned}
\right.\qquad
\end{equation}
for a vector-field source term  $\Xi : \Omega \to \R^d$. Because
$$ 
\Div (\nabla^s U) = \Delta U + \nabla \Div U,
$$
we see that \eqref{eq:elliptic-korn} is nothing but a Lamé-type system with a kind of homogeneous Robin (or mixed) boundary condition.

We define the Hilbert spaces
$$
\VV_1 :=   \left\{ W : \Omega \to \R^d \mid W \in H^1(\Omega), \; W \cdot n(x) = 0 \text{ on } \partial \Omega \right\}
$$
and 
$$
\VV_0 :=   \left\{ W : \Omega \to \R^d \mid W \in H^1(\Omega), \; W \cdot n(x) = 0 \text{ on } \partial \Omega , \; P_\Omega \la \nabla^a W \ra = 0 \right\},
$$
where $P_\Omega$ denotes the orthogonal projection onto the set $\AA_\Omega = \{ A \in \MM^a_d(\R); \;  Ax \in \RR_\Omega \}$ of all skew-symmetric matrices giving rise to a centered infinitesimal rigid displacement field preserving $\Omega$ (see \eqref{eq:RROmega} for the definition of $\RR_\Omega$).
Both spaces are endowed with the~$H^1(\Omega)$ norm. We then denote 
$$
\VV_\alpha := 
\begin{cases}
\VV_1 \quad\text{if}\quad \alpha \not\equiv 0 \\
\VV_0 \quad\text{if}\quad \alpha \equiv 0.
\end{cases}
$$
We also define on $\VV_\alpha$ the bilinear form 
$$
A_\alpha(U,W) :=   \int_{\Omega} \nabla^s U : \nabla^s W  \, \dx 
+ \int_{\partial \Omega} \frac{\alpha(x)}{2-\alpha(x)} U \cdot W \, \d\sigma_{\! x},
$$
where $M : N := \sum_{ij} m_{ij} n_{ij}$ for two matrices $M = (m_{ij})$, $N = (n_{ij})$.

\smallskip

The coercivity of the bilinear form $A_\alpha$ is related to Korn-type inequalities that we present below. 
We start stating  a first classical version of Korn's inequality:

\begin{lem}\label{lem:Korn1}
For any vector-field $U \in H^1(\Omega)$, we have
\begin{equation}\label{Korn1}
\inf_{R \in \RR} \| \nabla (U - R) \|^2 \lesssim \| \nabla^s U \|^2,  
\end{equation}
where we recall that $\RR$ is the space of all infinitesimal rigid displacement fields defined in~\eqref{eq:RR}, or equivalently, we have 
\begin{equation}\label{Korn1bis}
\| \nabla U \|^2 \lesssim \| \nabla^s U \|^2 +   |\langle \nabla^a U \rangle |^2.
\end{equation}
\end{lem}
For the statement of \eqref{Korn1} and its proof, we refer to \cite[Eq.~(1)]{DV02} where Friedrichs~\cite[Eq.~(13), Second case]{MR22750} and Duvaut-Lions~\cite[Eq.~(3.49)]{MR0521262} are quoted, as well as \cite[Theorem~2.2]{MR2119999} and the references therein. 

In the following lemma, we prove an estimate on $|\langle \nabla^a U \rangle |$ in the case $\alpha\not\equiv 0$.

\begin{lem}\label{lem:Korn2}
Supposing $\alpha\not\equiv 0$,  we have 
\begin{equation}\label{Korn1ter}
|\langle \nabla^a U \rangle |^2 \lesssim  \| \nabla^s U \|^2 +   \left\| \sqrt{\alpha \over {2-\alpha}} U \right\|_{L^2(\partial\Omega)}^2,
\end{equation}
for any vector-field $U \in H^1(\Omega)$. 
\end{lem}

\begin{proof}[Proof of Lemma~\ref{lem:Korn2}]
In order to establish \eqref{Korn1ter}, we argue by contradiction. We assume thus that  \eqref{Korn1ter} is not true, so that there exists a sequence $(U_n)_{n \in \N}$ in $H^1(\Omega)$ satisfying
$$
1 =  |\langle \nabla^a U_n \rangle |^2 \ge n \left(  \| \nabla^s U_n \|^2 +   \left\| \sqrt{\alpha \over {2-\alpha}} U_n \right\|^2_{L^2(\partial\Omega)} \right).
$$
Together with \eqref{Korn1bis} and \eqref{eq:PTypeineq} applied to each component of $U_n$, we obtain that $(U_n)_{n \in \N}$ is bounded in $H^1(\Omega)$. As a consequence, up to the extraction of a subsequence, there exists~$U \in H^1(\Omega)$ such that $U_n \wto U$ weakly in $H^1(\Omega)$ { and $U_n \to U$ strongly in $L^2(\Omega)$}. Passing to the limit in the above estimates satisfied by $(U_n)_{n \in \N}$, we get $|\langle \nabla^a U \rangle|^2 = 1$, $ \|  \sqrt{\alpha/(2-\alpha)} U \|^2_{L^2(\partial\Omega)} = 0$
and~$ \| \nabla^s U \| = 0$. From $\nabla^s U = 0$, we first deduce that there exist an antisymmetric matrix $A$ and a constant vector $b \in \R^d$ such that $U(x) = A x+b$ on $\Omega$,
and, thanks to the estimate $ \| \sqrt{\alpha/(2-\alpha)} U \|^2_{L^2(\partial\Omega)} = 0$, we deduce that 
$$
A x + b = 0 \quad \text{on} \quad \Gamma := \{ x \in \partial\Omega, \, \alpha(x) > 0 \},
$$ 
which has positive measure $|\Gamma|>0$ { using that $\alpha$ is a Lipschitz function}.  We fix $\bar x$ an interior point of $\Gamma$. As in the fourth step of the proof of Theorem~\ref{theo:Poisson}, we consider a family of smooth vector fields $a^1, \dots, a^d$ such that it is an orthonormal basis of $\R^d$ and  such that for any $x \in \partial\Omega$, $a^1(x)=n(x)$. We then introduce the flow $(\Phi^i_t)_{t \ge 0}$ associated to $a^i$ for $i= 2, \dots, d$. For $t$ small enough, $\Phi^i_t(\bar x)$ is still in the interior of $\Gamma$ so that 
$$
A a^i (\bar x) = {\d \over \d t} (A \Phi^i_t(\bar x) + b) = 0.
$$ 
Therefore, for any $i \ge 2$, one has, using that $A \bar x + b = 0$ so that $b = - A \bar x$
$$ 
a^i (\bar x) \cdot U(x) = a^i (\bar x) \cdot (Ax+b) = - A a^i (\bar x) \cdot x + A a^i (\bar x) \cdot \bar x =0,
$$ 
for any $x \in \Omega$, or, in other words, $U(x) \in \R  \bar n$ for any $x \in \Omega$, with $\bar n := n(\bar x)$. 
We may thus write $U(x) = \phi(x) \bar n$, with $\phi : \Omega \to \R$ an affine function, so that $\phi(x) = k \cdot x + k_0$, $k \in \R^d$, $k_0 \in \R$.
There exists next at least one index $i_0 \in \{ 1, \dots, d \}$ such that $\bar n_{i_0} \not = 0$ because $| \bar n | = 1$. Using again the fact that $\nabla^s U =0$ on~$\Omega$ and observing that $(\nabla U)_{ij} = k_i \bar n_j$, 
we deduce first $k_{i_0} = 0$ because $k_{i_0} \bar n_{i_0} = (\nabla^s U)_{i_0i_0} = 0$ and next $k_i = 0$ for any $i \not=i_0$ because $k_{i} \bar n_{i_0} =  2 (\nabla^s U)_{i_0 i} = 0$. We have thus established that  $U = n_0 := k_0 \, \bar n$ on $\Omega$, for some constant $n_0 \in \R^d$. We may alternatively prove that $\nabla U = 0$ and $U$ is constant again by using just the claim \cite[Eq.~(3)]{DV02}. Anyway, both arguments lead to the fact that $U=0$ because of the boundary condition on $\Gamma$ which is in contradiction with $|\langle \nabla^a U \rangle|^2 = 1$. That ends the proof of  \eqref{Korn1ter}. 
\end{proof}
Gathering \eqref{Korn1bis} and \eqref{Korn1ter}, we then have established the (probably classical) following Korn-type inequality:

\begin{lem}\label{lem:Korn3}
Assuming $\alpha \not\equiv 0$. For any vector-field $U \in H^1(\Omega)$, there holds
\begin{equation}\label{Korn1quar}
\| \nabla U \|^2 \lesssim  
\| \nabla^s U \|^2 +   \left\|  \sqrt{\alpha \over {2-\alpha}} U \right\|_{L^2(\partial\Omega)}^2.
\end{equation}
\end{lem}

For later reference, we also mention that a similar argument (and even a bit simpler, see also \cite[Eq.~(2)]{DV02} and \cite[Theorem~2.1]{MR2119999}) leads to the following variant of Korn's inequality: 
\begin{lem}\label{lem:Korn4}
For any vector-field $U \in H^1(\Omega)$, there holds
\begin{equation}\label{Korn1six}
 \| \nabla U \|^2 \lesssim  
  \| \nabla^s U \|^2 +   \|  U \|^2.
\end{equation}
\end{lem}

It is worth emphasizing that we also have the following Poincaré inequality: 
\begin{lem}\label{lem:Korn5}
For any $U \in H^1(\Omega)$ such that $U(x) \cdot n(x) = 0$ on $\partial \Omega$, there holds 
\begin{equation}\label{PoincareVectV0}
 \| U \|^2  \lesssim     \| \nabla U \|^2.
\end{equation}
\end{lem}

\begin{proof}[Proof of Lemma~\ref{lem:Korn5}]
 As before, we may argue by contradiction, assuming that  \eqref{PoincareVectV0} is not true, so that 
there exists a sequence $(U_n)_{n \in \N}$ in $H^1(\Omega)$ satisfying $U_n \cdot n(x) = 0$ on $\partial \Omega$ and such that
$$
1 =  \|  U_n \|^2  \ge n   \| \nabla U_n \|^2 .
$$
We immediately deduce that there exists $U \in H^1(\Omega)$ such that $\nabla U = 0$, $ \|  U  \|^2 = 1$ and $U \cdot n(x) = 0$ which gives our contradiction. 
\end{proof}

Gathering  \eqref{Korn1quar}  and  \eqref{PoincareVectV0}, we may state a last version of our first Korn inequality:
\begin{prop} 
Suppose that $\alpha \not\equiv 0$. For any $U \in H^1(\Omega)$ such that $U(x) \cdot n(x) = 0$ on~$\partial \Omega$, there holds 
\begin{equation}\label{PoincareKorn1}
 \| U \|_{H^1(\Omega)}^2  \lesssim    \| \nabla^s U \|^2 +   \left\|  \sqrt{\alpha \over {2-\alpha}} U \right\|_{L^2(\partial\Omega)}^2.
\end{equation}
\end{prop}

\smallskip

On the other hand, a less classical Korn's inequality has been established by Desvillettes and Villani~\cite{DV02}:
\begin{lem}\label{lem:Korn6}
For any vector-field $U \in H^1(\Omega)$ verifying $U \cdot n(x) = 0$ on $\partial \Omega$, one has 
\begin{equation}\label{Korn2}
\inf_{R \in \RR_\Omega} \| \nabla (U - R) \|^2 \lesssim \| \nabla^s U \|^2,
\end{equation}
where we remind that $\RR_\Omega$ stands for the space of centered infinitesimal rigid displacement fields defined in \eqref{eq:RROmega}, or equivalently one has
\begin{equation}\label{Korn2bis}
\| \nabla U \|^2 \lesssim  \| \nabla^s U \|^2 +   |P_\Omega \langle \nabla^a U \rangle |^2, 
\end{equation}
where we recall that $P_\Omega$ stands for the orthogonal projection onto the space $\AA_\Omega$ as defined before. 
\end{lem}
In the case when $\RR_\Omega = \{ 0 \}$, that is when $\Omega$ has no axi-symmetry, \eqref{Korn2} is nothing but the inequality stated in \cite[Theorem~3]{DV02} and for which a detailed constructive proof is provided therein. 
The proof of \eqref{Korn2} in the three dimensional case is also alluded in \cite[Section~5]{DV02}. We do not explain how the analysis developed in~\cite{DV02} makes possible  to get a constructive proof of \eqref{Korn2} in the general case (whatever is the dimension $d$), but rather briefly explain how \eqref{Korn2bis} may be established thanks to a compactness argument. 

\medskip
\noindent {\it Proof of~\eqref{Korn2bis}.} 
We first claim that for any vector-field $U \in H^1(\Omega)$ such that $U \cdot n(x) = 0$ on $\partial \Omega$, one has 
\begin{equation}\label{Korn2ter}
\| U \|^2 \lesssim  \| \nabla^s U \|^2 +   |P_\Omega \langle \nabla^a U \rangle |^2.
\end{equation}
Assume indeed by contradiction  that  \eqref{Korn2ter} is not true, so that 
there exists a sequence $(U_n)_{n \in \N}$ satisfying $U_n \cdot n(x) = 0$ on $\partial \Omega$ such that 
$$
1 =  \|  U_n \|^2 \ge n \left(  \| \nabla^s U_n \|^2 +  |P_\Omega \langle \nabla^a U_n \rangle |^2  \right).
$$
Together with the Korn inequality \eqref{Korn1six}, we deduce that there exists $U \in H^1(\Omega)$ satisfying $U \cdot n(x) = 0$ on $\partial \Omega$ such that (up to the extraction of a subsequence) $U_n \wto U$ weakly in~$H^1(\Omega)$ and $U_n \to U$ strongly in $L^2(\Omega)$. 
Passing to the limit in the estimates satisfied by~$(U_n)_{n \in \N}$, we first get $\nabla^s U = 0$ which implies that $U = Ax +b  \in \RR$. Moreover we obtain~$U \cdot n(x) = (Ax + b) \cdot n(x) = 0$ on $\partial \Omega$ and thus, thanks to the remark after \eqref{eq:RROmega} using the assumption \eqref{eq:Omega-centre}, we obtain that $b=0$ and hence $A \in \AA_\Omega$ or equivalently $Ax \in \RR_\Omega$. Finally, we also have $P_\Omega \la \nabla^a U \ra = P_\Omega A = 0$ which implies $A \in \AA_\Omega^\perp$ 
and thus~$A = 0$. We therefore obtain $U=0$ which is in contradiction with the fact that $\| U \|^2 = 1$. That ends the proof of \eqref{Korn2ter}. The proof of \eqref{Korn2bis} follows by gathering \eqref{Korn1six} and \eqref{Korn2ter}. 
\qed

\medskip

Gathering~\eqref{Korn2bis} with \eqref{Korn2ter}, we finally obtain the following Korn-type inequality: 
\begin{prop}
For any vector-field $U \in H^1(\Omega)$ such that $U \cdot n(x) = 0$ on $\partial \Omega$, there holds
\begin{equation}\label{PoincareKorn2}
\| U \|^2_{H^1(\Omega)} \lesssim  \| \nabla^s U \|^2 +   |P_\Omega \langle \nabla^a U \rangle |^2.
\end{equation}
\end{prop}

\smallskip

We can now state our result concerning the existence, uniqueness and regularity of solutions to the elliptic system \eqref{eq:elliptic-korn}.

\begin{theo}\label{theo:regH2-korn}
For any given $\Xi \in L^2(\Omega)$, there exists a unique solution $U \in \VV_\alpha$ to the variational problem associated to \eqref{eq:elliptic-korn}, namely
\begin{equation}\label{eq:elliptic-korn-var}
A_\alpha(U,W) = ( \Xi , W ) \quad\forall \, W \in \VV_\alpha.
\end{equation}
If furthermore $\Xi$ satisfies the condition $\la \Xi , Ax \ra = 0$ for any $Ax \in \RR_\Omega$ when $\alpha \equiv 0$, then the variational solution $U$ to \eqref{eq:elliptic-korn} satisfies $U \in H^2(\Omega)$ with 
$$
\| U \|_{H^2(\Omega)} \lesssim \| \Xi \|,
$$
and moreover $U$ verifies \eqref{eq:elliptic-korn} a.e.
\end{theo}

The proof of Theorem~\ref{theo:regH2-korn} follows the same steps as the proof of Theorem~\ref{theo:Poisson}. We briefly present it below.

\begin{proof}[Proof of Theorem~\ref{theo:regH2-korn}] 
We split the proof into four steps, the three last ones being devoted to the proof of the $H^2$ regularity estimate.

\smallskip\noindent{\sl Step 1. } 
Thanks to the above Korn-type inequalities, more precisely \eqref{PoincareKorn1} for the case $\alpha \not\equiv 0$ and \eqref{PoincareKorn2} for the case $\alpha \equiv 0$, we deduce that the bilinear form $A_\alpha$ is coercive in $\VV_\alpha$, that is, there is a constant $\lambda >0$ such that
$$
\forall \,  U \in \VV_\alpha, \quad  \quad  \lambda \left(\| U  \|^2 + \| \nabla U \|^2 \right) \le  A_\alpha(U,U).
$$
One can therefore apply Lax-Milgram theorem which gives us the existence and uniqueness of $U \in \VV_\alpha$ satisfying \eqref{eq:elliptic-korn-var}.

\smallskip

For the remainder of the proof, we additionally assume that $\la \Xi , Ax \ra = 0$ for any $Ax \in \RR_\Omega$ when $\alpha \equiv 0$. We then claim that \eqref{eq:elliptic-korn-var} can be improved into the following new variational formulation: there exists a unique $U \in \VV_\alpha$ verifying
\begin{equation}\label{eq:elliptic-korn-var-bis}
A_\alpha(U,W) = ( \Xi , W ), \quad\forall \, W \in \VV_1.
\end{equation}
In the case $\alpha \not \equiv 0$ or $\alpha \equiv 0$ with a non axi-symmetric domain $\Omega$, that is $\RR_\Omega = \{ 0 \}$, equation~\eqref{eq:elliptic-korn-var-bis} is nothing but \eqref{eq:elliptic-korn-var} since in these cases $\VV_\alpha = \VV_1$. When $\alpha \equiv 0$ and $\Omega$ has rotational symmetry, that is $\RR_\Omega \neq \{ 0 \}$, for any $W \in \VV_1$ we have $W - P_\Omega \la \nabla^a W \ra x \in \VV_0$ and therefore 
$$
\begin{aligned}
A_\alpha (U,W) 
&= A_\alpha (U , W - P_\Omega \la \nabla^a W \ra x) \\
&= \int_{\Omega} \Xi \cdot W \, \dx  
- \int_{\Omega} \Xi \cdot (P_\Omega \la \nabla^a W \ra x) \, \dx \\
&= \int_{\Omega} \Xi \cdot W \, \dx  
\end{aligned}
$$
where we have used that $\nabla^s (P_\Omega \la \nabla^a W \ra x) = 0$ in the first line, formulation \eqref{eq:elliptic-korn-var} in the second line, and the condition $\la \Xi , Ax \ra = 0$ for any $Ax \in \RR_\Omega$ in the third line, since~$P_\Omega \la \nabla^a W \ra x \in \RR_\Omega$ by definition.

\smallskip\noindent{\sl Step 2.} For any small enough open set $\omega \subset \Omega$, we fix 
a vector field $a \in C^2(\bar\Omega)$ such that $|a| = 1$ on $\omega$ and $a \cdot n = 0$ on $\partial \Omega$, 
and we set $X := a \cdot \nabla$ the associated differential operator. For a smooth solution $U$ to \eqref{eq:elliptic-korn-var-bis}, we compute
\bean
\| \nabla^s  X U\|^2  &=& 
 (\nabla^s  U, X^* \nabla^s  X U  ) +  ([\nabla^s,X]  U,   \nabla^s  X U  ) 
\\
&=& (\nabla^s  U,   \nabla^s  X^*X U  ) +  (\nabla^s  U, [X^*, \nabla^s ] X U  ) + ([\nabla^s,X]  U,   \nabla^s  X U  ) 
\eean
where we have used~\eqref{eq:X*}. 
On the other hand, we have the following formal equality
$$
\int_{\partial\Omega} (XU)\cdot (XU) \, \frac{\alpha}{2-\alpha}  \d\sigma_{\! x}
=  \int_{\partial\Omega}  \frac{\alpha}{2-\alpha} \, U \cdot ( X^* XU) \, \d\sigma_{\! x} 
- \int_{\partial\Omega}   \biggl( X \frac{\alpha}{2-\alpha} \biggr) \, U \cdot (XU) \, \d\sigma_{\! x} .
$$
We define 
\bean
(\AA W)_{ij} &:=& \frac12 \bigl( {[\partial_i,X]} W_j + {[\partial_j,X]} W_i \bigr)
\\
(\BB W)_{ij} &:=& \frac12 \bigl( {[{X^*}, {\partial_i}]}  W_j +{[{X^*}, {\partial_j}]}  W_i \bigr).
\eean
Supposing the additional regularity assumption $X^*XU \in \VV_1$, using { $(\nabla^s)^* \nabla^s = - \operatorname{div}(\nabla^s \cdot)$} and making the choice $W := X^*XU$ in the variational equation~\eqref{eq:elliptic-korn-var}, we obtain 
\begin{align*}
&\| \nabla^s  X U\|^2 + \int_{\partial\Omega} (XU)\cdot (XU) \, \frac{\alpha}{2-\alpha}  \d\sigma_{\! x} 
\\
&\quad=(\Xi, X^*X U  ) +  (\nabla^s  U, \BB X U  ) + (\AA  U,   \nabla^s  X U  ) -  \int_{\partial\Omega}   \biggl( X \frac{\alpha}{2-\alpha} \biggr) \, U \cdot (XU) \, \d\sigma_{\! x} .
\end{align*}
From the Korn inequalities \eqref{Korn1quar} (when $\alpha\not\equiv0$) and  \eqref{Korn1six} (when $\alpha\equiv0$), we first deduce
\bean
\| \nabla X U\|^2  &\lesssim  \| \Xi \| \| X^*X U \|  + \| \nabla  U \|  \| \BB X U \| + \| \AA  U \| \|  \nabla  X U  \|  \\
&\qquad \qquad + \|U\|_{L^2(\partial \Omega)} \|XU\|_{L^2(\partial\Omega)} + \|XU\|^2.
\eean

Then, since 
\bean
{[\partial_i,X]} = (\partial_i a) \cdot \nabla, \quad
{[{X^*}, {\partial_i}]} =  \partial_i (\hbox{div} a) + (\partial_i a) \cdot \nabla,
\eean
we deduce that 
$$
\| \AA W \|  +  \| \BB W \|   \lesssim \| W \|_{H^1(\Omega)}, \quad \forall \, W \in \VV_1.
$$
We also have the elementary estimates 
$$
\| X^*W \| +  \| X W \|   \lesssim  \| W \|_{H^1(\Omega)}, \quad \forall \, W \in \VV_1.
$$
Thanks to the already established estimate $\| U \|_{H^1(\Omega)} \lesssim \| \Xi \|$, we are then able to deduce that
\bean
\| \nabla X U\|^2  \lesssim \| \Xi \| \| \nabla X U\| + \| \Xi\|^2, 
\eean
and finally
\beqn\label{eq:RegDirX}
\| \nabla X U\|   \lesssim  \| \Xi \| .
\eeqn
Note that as in the proof of Theorem~\ref{theo:Poisson}, the multiplicative constants involved in our estimates depend on $\|a\|_{W^{2,\infty}(\Omega)}$ and $\|\alpha\|_{W^{1,\infty}(\Omega)}$.


\smallskip\noindent{\sl Step 3.} 
When we do not deal with an a priori smooth solution, but just with a solution $U \in \VV_\alpha$ to \eqref{eq:elliptic-korn-var-bis}, we modify the argument in the following way. 

We consider a small enough open set $\omega \in \Omega$, so that we may fix $a^1,\dots,a^d$ a family of smooth vector fields such that $(a^1,\dots,a^d)$ is an orthonormal basis of $\R^d$ at any point $x \in \omega$ and $a^1(x) = n(x)$ for any $x \in \partial \Omega \cap \partial \omega$. The construction of such a family is given in the Step 4 of the proof of Theorem \ref{theo:Poisson}. We set $A = (a^1,\dots,a^d)$. Let $k \in \{2, \dots, d\}$. Then $a = a^k$ is as in Step 2 and we define $\Phi_t$ the associated flow introduced in~\eqref{eq:flow}.

We define $J^h (x) := A(\Phi_h(x)) A(x)^{-1}$, so that in particular 
$J^h (x) n (x)  = n(\Phi_h(x))$ for any $h$. We next define 
$$
X^hU(x) := \frac{1}{h} \left( {^T\!}J^h(x) U(\Phi_h(x)) - U(x)\right),
$$
so that $X^h U \in \VV_1$ if $U \in \VV_\alpha$. Repeating the argument of Step 2, we get
\bean
\| \nabla^s  X^h  U\|^2  = (\nabla^s U, \nabla^s X^{h*}X^h U  ) +  (\nabla^s  U, \BB^h X^h U  ) + (\AA ^h  U,   \nabla^s  X^h  U  ), 
\eean
where we denote
\bean
X^{h*} M (x) &:=& {1 \over h} [ |\operatorname{det} D\Phi_{-h}(x)| J^h(\Phi_{-h}(x)) M(\Phi_{-h}(x)) - M(x)]
\\
(\AA^h  W)_{ij} &:=& \frac12 \bigl( {[\partial_i,X^h]} W_j + {[\partial_j,X^h]} W_i \bigr)
\\
(\BB^h  W)_{ij} &:=& \frac12 \bigl( {[{X^{h*}}, {\partial_i}]}  W_j +{[{X^{h*}}, {\partial_j}]}  W_i \bigr).
\eean
On the other hand, we have
\begin{align*}
&\int_{\partial \Omega} \frac{\alpha}{2-\alpha}(x) U(x) \cdot X^{h*} X^h U(x) \d\sigma_{\! x} \\
&\quad = \int_{\partial \Omega} \frac{\alpha}{2 - \alpha}(\Phi_h(x)) (X^hU)(x) \cdot (X^hU)(x) \d\sigma_{\! x} + \int_{\partial \Omega} U(x) \cdot X^hU(x) Y^h\bigg(\frac{\alpha}{2-\alpha}\bigg)(x) \d\sigma_{\! x},
\end{align*}
where 
\[ Y^h M(x) := \frac{1}{h} \Big(M(\Phi_h(x)) - M(x)\Big). \]
We also have that if $U \in \VV_\alpha$ then $X^{h*}X^h U\in \VV_1$ too. Indeed, we compute
\begin{align*}
X^{h*}X^hU (x)&= {1 \over h^2} |\operatorname{det} D\Phi_{-h}(x)| J^h(\Phi_{-h}(x)) \left(\left({^T\!J^h}(\Phi_{-h}(x))\right) U(x) - U(\Phi_{-h}(x))\right) \\
&\quad - {1 \over h^2} \left( {^T\!J^h}(x) U(\Phi_h(x)) - U(x)\right) =: T_1(x) + T_2(x),
\end{align*}
the last equality standing for a definition of $T_1$ and $T_2$. As already noticed, if $U \in \VV_\alpha$, then $X^hU(x) \cdot n(x)=0$ so that $T_2(x) \cdot n(x) =0$. 
Concerning~$T_1$, we first have 
$$
J^h(\Phi_{-h}(x)) \left({^T\!J^h}(\Phi_{-h}(x))\right) U(x) \cdot n(x) = U(x) \cdot n(x) = 0.
$$
Then, we remark that $J^h(\Phi_{-h}(x)) = {^T\!}J^{-h}(x)$, so that 
$$
J^h(\Phi_{-h}(x))  U(\Phi_{-h}(x)) \cdot n(x) = U(\Phi_{-h}(x)) \cdot J^{-h}(x) n(x) = U(\Phi_{-h}(x)) \cdot n(\Phi_{-h}(x))=0.
$$
Using this and the fact that $U$ is a solution of \eqref{eq:elliptic-korn-var-bis}, we deduce that
\begin{align*}
&\| \nabla^s  X^h  U\|^2 + \int_{\partial \Omega} \frac{\alpha}{2 - \alpha}(\Phi_h(x)) (X^hU)(x) \cdot (X^hU)(x) \d\sigma_{\! x} \\
&\quad = (\Xi, X^{h*}X^h U  ) +  (\nabla^s  U, \BB^h X^h U  ) + (\AA ^h  U,   \nabla^s  X^h  U  ) - \int_{\partial \Omega} U \cdot (X^h U) \, \bigg( Y^h \frac{\alpha}{2 - \alpha}\bigg) \d\sigma_{\! x}.
\end{align*}
Similarly as in the proof of Theorem~\ref{theo:Poisson}, one can prove the following elementary estimate
$$
\|X^h W\| + \| X^{h*}W \| +  \| \AA^h  W \|  +  \| \BB^h W \|   \lesssim \|W \|_{H^1(\Omega)}, \quad \forall \, W \in \VV_1.
$$
Using these bounds combined with the already established estimate $ \| U \|_{H^1(\Omega)} \lesssim \| \Xi \|$ and the Korn inequality, we deduce, as in the Poisson case, that
$$
\| \nabla X^h U\|   \lesssim  \| \Xi \|, \quad \forall \, |h| \le 1. 
$$ 
Passing to the limit $h \to 0$, we then get 
$$
\| \nabla X^0 U\|   \lesssim  \| \Xi \|, 
$$
with $X^0 U_j = a \cdot \nabla U_j + A \, (a \cdot \nabla A^{-1}) U_j$ for $j=1,\dots,d$. Note that as in the Poisson case, the multiplicative constants are uniform in $|h| \le 1$ and depend on $\|a\|_{W^{2,\infty}}$ and~$\|\alpha\|_{W^{1,\infty}}$. We then recover~\eqref{eq:RegDirX} by observing that we have $ \| A \, (a \cdot \nabla A^{-1}) U \|_{H^1} \lesssim  \| \Xi \|$.

\medskip\noindent
{\sl Step 4. } We set now $X_i := a^i \cdot \nabla$. From the second step, we have 
\beqn\label{eq:RegDirXi}
\| \nabla X_i U\|   \lesssim  \| \Xi \|, \quad \forall \, i = 2, \dots, d.
\eeqn
We first notice that 
$$
\partial_j = \sum_i a^i_jX_i = -\sum_i X_i^*(a^i_j \cdot).
$$
Combining this with~\eqref{eq:Deltau}, we deduce that 
\begin{align*}
\Xi_j &= -\Delta U_j - \partial_j(\operatorname{div} U) = \sum_i X_i^*X_i U_j + \sum_{i,\ell,m} X_i^*(a^i_ja^m_\ell X_m U_\ell) \\
&= X_1^* X_1 U_j + \sum_\ell X_1^*(a^1_ja^1_\ell X_1 U_\ell) + \sum_{i \neq 1} X_i^*X_i U_j+ \sum_{(i,m) \neq (1,1)} \sum_\ell X_i^*(a^i_ja^m_\ell X_m U_\ell). 
\end{align*}
We notice that $X_i^*(fg) = (X_i^*f) g - f(X_i g)$. Using then~\eqref{eq:RegDirXi} combined with the fact that for $i=1,\dots,d$, we have $a^i \in W^{2,\infty}(\Omega)$, 
we deduce
\beqn \label{eq:X1*X1_1}
X_1^*X_1 U_j + \sum_\ell a^1_j a^1_\ell X_1^*X_1 U_\ell = R_j(U,\Xi) \quad \text{with} \quad \|R_j(U,\Xi)\| \lesssim \|\Xi\|. 
\eeqn
Multiplying the equality in~\eqref{eq:X1*X1_1} by $a_j^1$ and then summing it over $j$, we get 
$$
2 \sum_\ell a^1_\ell X_1^*X_1 U_\ell = \sum_j a^1_j R_j(U,\Xi), 
$$
and thus 
\beqn \label{eq:a1X1*X1}
\big\|a^1 \cdot X_1^*X_1 U\big\| \lesssim \|\Xi\|. 
\eeqn
Coming back to~\eqref{eq:X1*X1_1} and using once more that $\delta_{j\ell} = a_j \cdot a_\ell$, so that
\beqn \label{eq:X1*X1}
 X_1^* X_1 U_j = \sum_{\ell,m} a^m_j a^m_\ell X_1^* X_1 U_\ell, 
 \eeqn
 we obtain that 
$$
\sum_{ m \neq 1, \ell \in \{1,\dots,d\}} a^m_ja^m_\ell X_1^*X_1U_\ell = R_j(U,\Xi) - 2 \sum_\ell a^1_j a^1_\ell X_1^*X_1 U_\ell . 
$$ 
Together with~\eqref{eq:a1X1*X1} and the fact that $\|R_j(U,\Xi)\|\lesssim \|\Xi\|$, it yields
\beqn \label{eq:amX1*X1}
\bigg\|\sum_{\ell, m \neq 1} a^m_ja^m_\ell X_1^*X_1U_\ell \bigg\| \lesssim \|\Xi\|. 
\eeqn 
Finally, using again~\eqref{eq:X1*X1},~\eqref{eq:a1X1*X1} and~\eqref{eq:amX1*X1} imply
$$
\|X_1^*X_1 U_j\| \lesssim \|\Xi\|. 
$$
Recalling that $[X_1,X_1^*] u = (a^1 \cdot \nabla \hbox{div} (a^1))u $, because { $\|U\|_{H^1(\Omega)} \lesssim \|\Xi\|$}, the above inequality implies 
$$
\| X^2_1 U\|   \lesssim  \| \Xi \|,
$$
and then together with \eqref{eq:RegDirXi}, we have established
$$
\| X_i X_j U\|   \lesssim  \| \Xi \|, \quad \forall \, i ,j = 1, \dots, d.
$$
We can then conclude the proof of Theorem~\ref{theo:regH2-korn} as in the one of Theorem~\ref{theo:Poisson}. 
\end{proof}


\section{Proof of Theorem~\ref{theo:hypo}}\label{sec:proof}
\label{sec:proofHypo}

Consider the operator $\LLL$ defined in~\eqref{eq:dtf=Lf}. For any $f \in \HH $ we decompose $f = \pi f + f^\perp$ with the macroscopic part $\pi f$ given by
$$
\pi f (x,v) = \varrho(x) \mu(v) + m(x) \cdot v \mu(v) + \theta(x) \, \frac{(|v|^2 - d)}{\sqrt{2d}} \, \mu(v),
$$
where the mass, momentum and energy are defined respectively by
$$
\varrho(x) = \int_{\R^d} f (x,v) \, \dv , \quad 
m(x) = \int_{\R^d} v f (x,v) \, \dv  
\quad \text{and} \quad
\theta(x) = \int_{\R^d} \frac{(|v|^2 - d)}{\sqrt{2d}} \, f (x,v) \, \dv .
$$
Remark that 
$$
\| f \|_{\HH}^2 
= \| f^\perp \|_{\HH}^2 + 
\| \pi f \|_{\HH}^2
$$
and
$$ 
\| \pi f \|_{\HH}^2 
= \| \varrho \|_{L^2_x(\Omega)}^2 
+ \| m \|_{L^2_x(\Omega)}^2 
+ \| \theta \|_{L^2_x(\Omega)}^2.
$$

\smallskip

The focus of the remainder of this section will be the proof of Theorem~\ref{theo:hypo} (note that Theorem~\ref{theo:main} is a direct consequence of Theorem~\ref{theo:hypo}). As explained in Subsection~\ref{ssec:main}, in Theorem~\ref{theo:hypo}, the construction of the scalar product $\la\!\la \cdot , \cdot \ra \! \ra$ on the space $\HH$ begins with the usual scalar product, which gives us a control of the microscopic part $f^\perp$, and after that, step by step, new terms are added to it in order to control all components of the macroscopic part $\pi f$. The construction of each of those terms is performed from Section~\ref{ssec:micro} through Section~\ref{ssec:mass}, and then in Section~\ref{ssec:conclusion} we shall complete the proof of Theorem~\ref{theo:hypo}.

\smallskip
We consider hereafter $f$ satisfying the conditions of Theorem~\ref{theo:hypo}, namely 
$f \in \mathrm{Dom}(\LLL)$ satisfying the boundary condition \eqref{eq:BdyCond}, so that in particular \eqref{eq:invariantsBoundary1} holds, which translates into
\beqn\label{eq:mcdotn=0}
m (x) \cdot n(x) = 0 \quad \text{for} \quad x \in \partial \Omega,
\eeqn
and satisfying assumption~\eqref{eq:C1} which means
$$
\langle \varrho \rangle =  \int_{\Omega} \varrho \, \dx = 0.
$$
In the specular reflection case ($\alpha\equiv 0$ in \eqref{eq:BdyCond}), the additional assumptions \eqref{eq:C2}-\eqref{eq:C3} hold, which corresponds to
\beqn\label{eq:C2C3pif}
\langle \theta \rangle = \int_{\Omega} \theta \, \dx = 0
\quad\text{and}\quad 
\langle R \cdot m \rangle = \int_{\Omega} R \cdot m \, \dx =  0 \quad \forall \, R \in \RR_\Omega.
\eeqn

For simplicity we introduce the notations $f_{\pm} := \gamma_{\pm} f$, $D^\perp := \mathrm{Id} - D$, where $D$ is given by~\eqref{eq:def_D} and $\partial \HH_+ := L^2(\Sigma_+ ; \mu^{-1}(v) n(x) \cdot v)$. 
It is worth emphasizing that because $f \in \mathrm{Dom}(\LLL)$, the trace functions $f_{\pm} $ are well defined. We refer the interested reader to  \cite{MR274925,MR777741} for the classical definition of the trace of a solution to a transport equation as well as to \cite{MR1765137,MischlerCMP2000,MR2150445} for a more modern approach.

\subsection{Microscopic part}\label{ssec:micro}

We start with the following result, giving a control of the microscopic part $f^\perp$ and a boundary term.

\begin{lem}\label{lem:micro}
There exists $\lambda >0$ such that 
$$
\la - \LLL f , f \ra_{\HH} \ge \lambda \| f^\perp \|_{\HH}^2 + \frac12 \| \sqrt{\alpha(2-\alpha)} D^\perp f_+ \|_{\partial \HH_+}^2.
$$   
\end{lem}

\begin{proof}[Proof of Lemma~\ref{lem:micro}]
We write 
$$
\la - \LLL f , f \ra_{\HH} = \la - \CCC f , f \ra_{\HH} 
+ \la  v \cdot \nabla_x f , f \ra_{\HH}.
$$    
Thanks to \eqref{eq:coercivity} one has 
$$
\la - \CCC f , f \ra_{\HH} \ge  \lambda \| f^\perp \|_{\HH}^2.
$$

For the second term, we first get thanks to an integration by parts
$$
\la  v \cdot \nabla_x f , f \ra_{\HH}  
= \int_\OO (v \cdot \nabla_x f) f \mu^{-1} \, \dx \, \dv 
= \frac12 \int_{\Sigma} \gamma f^2 \mu^{-1} n(x) \cdot v \, \d\sigma_{\! x} \, \dv.
$$
Writing $\gamma f^2 =  f_{+}^2 \mathbf 1_{\Sigma_{+}} + f_{-}^2 \mathbf 1_{\Sigma_{-}}$ and using the boundary condition \eqref{eq:BdyCond}, we thus obtain
$$
\begin{aligned}
\la  v \cdot \nabla_x f , f \ra_{\HH}  
&= \frac12 \int_{\Sigma_{+}}  f_{+}^2 \mu^{-1} |n(x) \cdot v| \, \d\sigma_{\! x} \, \dv 
- \frac12 \int_{\Sigma_{-}}  f_{-}^2 \mu^{-1} |n(x) \cdot v| \, \d\sigma_{\! x} \, \dv \\
&= \frac12 \int_{\Sigma_{+}}  f_{+}^2 \mu^{-1} |n(x) \cdot v |\, \d\sigma_{\! x} \, \dv \\
&\quad 
- \frac12 \int_{\Sigma_{-}}  \big \{ (1-\alpha(x))f_{+}(x,R_x v) + \alpha(x) D f_{+}(x,v)  \big\}^2 \mu^{-1} |n(x) \cdot v| \, \d\sigma_{\! x} \, \dv .
\end{aligned}
$$
We apply the change of variables $v \mapsto R_x v$, so that $\Sigma_{-}$ transforms into $\Sigma_{+}$, which yields
$$
\begin{aligned}
\la  v \cdot \nabla_x f , f \ra_{\HH}  
&= \frac12 \int_{\Sigma_{+}}  f_{+}^2 \mu^{-1} |n(x) \cdot v| \, \d\sigma_{\! x} \, \dv \\
&\quad
- \frac12 \int_{\Sigma_{+}} \big \{ (1-\alpha(x))f_{+} + \alpha(x) D f_{+}  \big \}^2 \mu^{-1} |n(x) \cdot v |\, \d\sigma_{\! x} \, \dv,
\end{aligned}
$$
since $Df_{+}(x,R_x v) = Df_{+}(x,v)$ and $|n(x) \cdot R_x v| = |n(x) \cdot v|$.
Writing $f_{+} = D^\perp f_{+} + D f_{+}$, one has 
$$
\int_{\Sigma_+} f_+^2 \mu^{-1}  n(x) \cdot v \, \d\sigma_{\! x} \, \dv 
= \int_{\Sigma_+} (Df_+)^2 \mu^{-1}  n(x) \cdot v \, \d\sigma_{\! x} \, \dv 
+\int_{\Sigma_+} (D^\perp f_+)^2 \mu^{-1}  n(x) \cdot v \, \d\sigma_{\! x} \, \dv,
$$ 
since $Df_+ \perp D^\perp f_+ $ in $\partial\HH_+$. 
All together, we conclude to
$$
\begin{aligned}
&\la  v \cdot \nabla_x f , f \ra_{\HH}  \\
&\qquad 
= \frac12 \int_{\Sigma_{+}}  \left\{ (Df_{+})^2 +  (D^\perp f_{+})^2 - [ (1-\alpha(x)) D^\perp f_{+} +  D f_{+} ]^2\right\} \mu^{-1} n(x) \cdot v \, \d\sigma_{\! x} \, \dv \\
&\qquad= \frac12 \int_{\Sigma_{+}}  \left\{ [1- (1-\alpha(x))^2] (D^\perp f_{+})^2 - 2(1-\alpha(x))  Df_+ D^\perp f_+  \right\} \mu^{-1} n(x) \cdot v \, \d\sigma_{\! x} \, \dv \\
&\qquad= \frac12 \int_{\Sigma_{+}} \alpha(x)(2-\alpha(x)) (D^\perp f_{+})^2\mu^{-1} n(x) \cdot v \, \d\sigma_{\! x} \, \dv.
\end{aligned}
$$
We finish the proof by gathering previous estimates.
\end{proof}

\subsection{Boundary terms}\label{ssec:bdry}

We start by stating a technical lemma which will be useful to treat the boundary terms in what follows. 

\begin{lem}\label{lem:boundary}
Let  $\phi: \R^d \to \R$. For any $x \in \partial\Omega$, there holds
$$
\begin{aligned}
\int_{\R^d}  \phi(v) \gamma f(x,v) \, n(x) \cdot v \, \dv  
& = \int_{\Sigma^x_{+}}  \phi(v) \alpha(x) D^\perp f_{+} \, n(x) \cdot v \, \dv  \\
&\quad 
+ \int_{\Sigma^x_{+}}  \left\{ \phi(v) - \phi(R_x v) \right\}  (1-\alpha(x)) D^\perp f_{+} \, n(x) \cdot v \, \dv  \\
&\quad 
+ \int_{\Sigma^x_{+}}   \left\{ \phi(v) - \phi(R_x v) \right\}  D f_{+} \, n(x) \cdot v \, \dv .
\end{aligned}
$$    
\end{lem}

\begin{proof}[Proof of Lemma~\ref{lem:boundary}]
We first write, thanks to the decomposition $\gamma f =  f_{+} \mathbf 1_{\Sigma_{+}} + f_{-} \mathbf 1_{\Sigma_{-}}$, 
$$
\begin{aligned}
\int_{\R^d} \phi(v) \gamma f(x,v) \, n(x) \cdot v \, \dv 
= \int_{\Sigma^x_{+}}   \phi(v) f_{+} \, n(x) \cdot v \, \dv 
- \int_{\Sigma^x_{-}}   \phi(v) f_{-} \, |n(x) \cdot v| \, \dv .
\end{aligned}
$$       
Applying the boundary condition \eqref{eq:BdyCond} and then the change of variables $v \mapsto R_x v$, we hence obtain 
$$
\begin{aligned}
&\int_{\Sigma^x_{-}}  \phi(v) f_{-} \, |n(x) \cdot v| \, \dv\\
&\qquad 
= \int_{\Sigma^x_{-}}  \phi(v)  \left\{ (1-\alpha(x)) f_{+}(x,R_x v) + \alpha(x) D f_{+} (x,v) \right \} \, |n(x) \cdot v| \, \dv \\
&\qquad 
= \int_{\Sigma^x_{+}}  \phi(R_x v) \left\{ (1-\alpha(x)) f_{+}(x,v) + \alpha(x) D f_{+} (x,v)  \right\} \, |n(x) \cdot v| \, \dv  ,
\end{aligned}
$$
since $Df_{+}(x,R_x v) = Df_{+}(x,v)$ and $|n(x) \cdot R_x v| = |n(x) \cdot v|$. We write $f_{+} = D^\perp f_{+} + D f_{+} $ and thus
$$ 
\begin{aligned}
&\int_{\R^d} \phi(v) \gamma f(x,v) \, n(x) \cdot v \, \dv \\
&\qquad 
= \int_{\Sigma^x_{+}} 
\left\{ \phi(v) f_{+} - \phi(R_x v)(1-\alpha(x)) f_{+} - \phi(R_x v) \alpha(x) D f_{+}  \right\}  n(x) \cdot v \, \dv \\
&\qquad 
= \int_{\Sigma^x_{+}} \phi(v) \alpha(x) D^\perp f_{+} \, n(x) \cdot v \, \dv \\
&\qquad\quad 
+ \int_{\Sigma^x_{+}} \left\{ \phi(v) - \phi(R_x v) \right\}  (1-\alpha(x)) D^\perp f_{+} \, n(x) \cdot v \, \dv \\
&\qquad\quad 
+ \int_{\Sigma^x_{+}} \left\{ \phi(v) - \phi(R_x v) \right\}  D f_{+} \, n(x) \cdot v \, \dv,
\end{aligned}
$$ 
which concludes the proof.
\end{proof}

\subsection{Energy}\label{ssec:energy}

In this subsection we construct a functional in order to control the energy component of the macroscopic part $\pi f$.
We denote 
$$
\theta[g] := \int_{\R^d} \frac{(|v|^2 - d)}{\sqrt{2d}} \, g\, \dv,
$$
so that $\theta = \theta[f]$. 
We define $u[\theta]$ as the solution to the elliptic equation \eqref{eq:elliptic} associated to $\xi=\theta \in L^2_x(\Omega)$ given by Theorem~\ref{theo:Poisson}, in particular
\beqn \label{eq:uthetafH2}
\| u[\theta] \|_{H^2_x(\Omega)} \lesssim \| \theta \|_{L^2_x(\Omega)}.
\eeqn 
It is worth noticing that in the specular reflection case, that is when $\alpha \equiv 0$ in \eqref{eq:BdyCond}, we have $\la \theta \ra = 0$ from \eqref{eq:C2C3pif}, so that the solution $u[\theta]$ to the Poisson equation with Neumann boundary condition is well-defined.
 
\smallskip
 
We also introduce the vector $p = (p_i)_{1 \le i \le d}$ defined by 
$$
p_i(v): = v_i \, \frac{(|v|^2 - d-2)}{\sqrt{2d}},
$$
and the associated moment functional  $M_p[g] = ( M_{p_i} [g])_{1 \le i \le d}$ given by 
\begin{equation}\label{eq:def-Mp}
M_{p_i}[g] = \int_{\R^d} v_i \, \frac{(|v|^2 - d-2)}{\sqrt{2d}} \, g \, \dv.
\end{equation}

\begin{lem}\label{lem:uthetaL} 
One has  
\begin{equation}\label{thetaLf}
\theta [\LLL f] 
= -\sqrt{\frac{2}{d}} \, \nabla_x \cdot m 
- \nabla_x \cdot M_p[f]
\end{equation}
and 
\begin{equation}\label{Mpf}
M_p[f] = M_p[f^\perp].
\end{equation}
As a consequence, from Theorem~\ref{theo:Poisson}, the unique variational solution $u[\theta [\LLL f]]$ to \eqref{eq:elliptic} associated to~$\xi=\theta[ \LLL f]$ satisfies 
\beqn \label{eq:uthetaLfH1}
\| u[\theta [\LLL f]] \|_{H^1_x(\Omega)}  
\lesssim \| m \|_{L^2_x(\Omega)} + \| f^\perp \|_{\HH} + \| \sqrt{\alpha(2-\alpha)} \, D^\perp f_+ \|_{\partial \HH_+}. 
\eeqn
\end{lem}

\begin{proof}[Proof of Lemma~\ref{lem:uthetaL}]
We start by proving \eqref{thetaLf}. By writing $\LLL f = - v \cdot \nabla_x f + \CCC f^\perp$ we have $\theta[\LLL f] = \theta [ - v \cdot \nabla_x f]$. We then compute
$$
\begin{aligned}
\theta[ - v \cdot \nabla_x f] 
&= -   \nabla_x \cdot \int_{\R^d} \frac{(|v|^2 - d)}{\sqrt{2d}} v f \, \dv  \\
&= -  \sqrt{\frac{2}{d}} \nabla_x \cdot \int_{\R^d} v f \, \dv
-  \nabla_x \cdot \int_{\R^d} \frac{(|v|^2 - d-2)}{\sqrt{2d}} v f \, \dv , 
\end{aligned}
$$
and this concludes the proof of \eqref{thetaLf}. Moreover, using the decomposition
\beqn\label{eq:fdecomposition}
f = \varrho \mu + m\cdot v \mu + \theta \frac{|v|^2-d}{\sqrt{2d}} \mu + f^\perp,
\eeqn
a straightforward computation  gives 
$$
\begin{aligned}
M_p[f] 
&= \varrho \int_{\R^d} p(v)  \mu \, \dv 
+ m_i \int_{\R^d} v_i p(v) \mu  \, \dv
+ \theta \int_{\R^d} p(v)  \left(\frac{|v|^2 - d}{\sqrt{2d}} \right) \mu  \, \dv
+ M_p[f^\perp] .
\end{aligned}
$$
We conclude to  \eqref{Mpf}, since $\int_{\R^d} p(v)  \mu \, \dv = \int_{\R^d} v_i p(v)  \mu \, \dv  = \int_{\R^d} (|v|^2-d)p(v)  \mu \, \dv = 0$.

\smallskip
From Theorem~\ref{theo:Poisson}, there exists a unique variational solution $u := u[\theta [\LLL f]]$ to \eqref{eq:elliptic} associated to~$\xi=\theta[ \LLL f]$. 
Thanks to Step~1 in the proof of Theorem~\ref{theo:Poisson}, this solution satisfies 
\begin{equation}\label{eq:uthetaLf1}
\lambda \| u \|_{H^1_x (\Omega)}^2
\le \| \nabla_x u \|_{L^2_x(\Omega)}^2
+ \| \sqrt{\tfrac{\alpha}{2-\alpha}} \, u \|_{L^2_x(\partial\Omega)}^2 , 
\end{equation}
for some constant $\lambda >0$. Moreover, thanks to the variational formulation \eqref{eq:varValpha}, one has
\bean
&& \| \nabla_x u \|_{L^2_x(\Omega)}^2
+ \| \sqrt{\tfrac{\alpha}{2-\alpha}} \, u \|_{L^2_x(\partial\Omega)}^2  \\
&&\qquad 
= -\int_{\Omega} \left(\sqrt{\frac{2}{d}} \, \nabla_x \cdot m 
+ \nabla_x \cdot M_p[f] \right) u \, \dx \\
&&\qquad 
=  \int_{\Omega} \left(\sqrt{\frac{2}{d}} \, m 
+  M_p[f] \right) \cdot \nabla_x  u \, \dx
-  \int_{\partial\Omega} \left(\sqrt{\frac{2}{d}} \, m 
+  M_p[f] \right) \cdot n(x) \,  u \, \d\sigma_{\! x},
\eean
where we have performed one integration by parts in the second equality. As a consequence, we have 
\begin{equation}\label{eq:uthetaLf2}
\begin{aligned}
& \| \nabla_x u \|_{L^2_x(\Omega)}^2
+ \| \sqrt{\tfrac{\alpha}{2-\alpha}} \, u \|_{L^2_x(\partial\Omega)}^2  \\
&\qquad 
= \int_{\Omega} \left(\sqrt{\frac{2}{d}} \, m 
+  M_p[f^\perp] \right) \cdot \nabla_x  u \, \dx
-  \int_{\partial\Omega} M_p[f] \cdot n(x) \,  u \, \d\sigma_{\! x}, 
\end{aligned}
\end{equation}
where we have used \eqref{Mpf} and that $m \cdot n = 0$ as noticed in \eqref{eq:mcdotn=0}.
For the boundary term appearing in last equation, we observe that thanks to Lemma~\ref{lem:boundary} and because $|v|^2 = |R_x v|^2$, for any $x \in \partial \Omega$, we have
$$
\begin{aligned} 
M_p[f] \cdot n(x)
&= \int_{\R^d} \frac{|v|^2-d-2}{\sqrt{2d}} f   \, n(x) \cdot v \, \dv \\
&=  \alpha(x) \int_{\Sigma_+^x} \frac{|v|^2-d-2}{\sqrt{2d}}\,  D^\perp f_+ \, n(x) \cdot v \, \dv,
\end{aligned}
$$
and therefore 
$$
\left| \int_{\partial\Omega} M_p[f] \cdot n(x) \,  u \, \d\sigma_{\! x} \right|
\lesssim  \Big\| \sqrt{\alpha(2-\alpha)} \, D^\perp f_+ \Big\|_{\partial \HH_+} \Big\| \sqrt{\tfrac{\alpha}{2-\alpha}} \, u \Big\|_{L^2_x(\partial\Omega)}.
$$
Remarking that
$$
\| M_p[f^\perp] \|_{L^2_x(\Omega)} \lesssim \| f^\perp \|_{\HH},
$$
we finally obtain \eqref{eq:uthetaLfH1} by gathering the above estimate on the boundary term together with \eqref{eq:uthetaLf1} and \eqref{eq:uthetaLf2}, and using Cauchy-Schwarz inequality.
\end{proof}

\smallskip

We next establish the following result, which gives us a control of the energy $\theta$.

\begin{lem}\label{lem:energy}
There are constants $\kappa_1, C >0$ such that 
$$
\begin{aligned}
&\la -\nabla_x u[\theta] ,  M_p [\LLL f] \ra_{L^2_x(\Omega)}
+ \la -\nabla_x u[\theta [\LLL f]],  M_p [f]  \ra_{L^2_x(\Omega)}\\
&\qquad\qquad
\ge \kappa_1 \| \theta \|_{L^2_x(\Omega)}^2 
- C \| m \|_{L^2_x(\Omega)} \| f^\perp \|_{\HH} 
- C \| f^\perp \|_{\HH}^2 
-C \| \sqrt{\alpha(2-\alpha)} D^\perp f_+ \|_{\partial \HH_+}^2.
\end{aligned}
$$
\end{lem}

\begin{proof}[Proof of Lemma~\ref{lem:energy}]
Using~\eqref{eq:uthetaLfH1} and~\eqref{Mpf}, one has
$$
\begin{aligned}
\left| \la -\nabla_x u[\theta [\LLL f]],  M_p [f^\perp]  \ra_{L^2_x(\Omega)} \right|
&\lesssim  \| \nabla_x u[\theta [\LLL f]]\|_{L^2_x(\Omega)} \| f^\perp \|_{\HH} \\
&\lesssim \| m \|_{L^2_x(\Omega)} \| f^\perp \|_{\HH} 
+ \| f^\perp \|_{\HH}^2 +  \| \sqrt{\alpha(2-\alpha)} D^\perp f_+ \|_{\partial \HH_+}^2,
\end{aligned}
$$
which allows us to bound the second term in the LHS of the estimate of the statement. 
For the first term, writing $M_p[\LLL f] = M_p[- v \cdot \nabla_x f] + M_p[ \CCC f^\perp]$ one obtains 
$$
\la -\nabla_x u[\theta] ,  M_p [\LLL f] \ra_{L^2_x(\Omega)} 
= T_1 + T_2 
$$
with 
$$
T_1 :=  \la \partial_{x_i} u[\theta] , \partial_{x_j} \int_{\R^d} p_i(v) v_j f \, \dv \ra_{L^2_x(\Omega)}
$$
and
$$
T_2 :=  \la - \nabla_x u[\theta] , \int_{\R^d} p(v) \CCC f^\perp \, \dv  \ra_{L^2_x(\Omega)}.
$$
For the term $T_2$, we remark that 
$$
\int_{\R^d} p(v) \CCC f^\perp \, \dv = \left( f^\perp , \CCC (p \mu) \right)_{L^2_v(\mu^{-1})},
$$
so that from the property (A3) on $\CCC$ and~\eqref{eq:uthetafH2}, we get
$$
|T_2| 
\lesssim \| \nabla_x u[\theta] \|_{L^2_x(\Omega)} \| f^\perp \|_{\HH} 
\lesssim \| \theta \|_{L^2_x(\Omega)} \| f^\perp \|_{\HH}.
$$
For the term $T_1$, we write
$$
\begin{aligned}
T_1 
&= - \la \partial_{x_j} \partial_{x_i} u[\theta] , \int_{\R^d} p_i(v) v_j f \, \dv \ra_{L^2_x(\Omega)} 
+ \int_{\partial\Omega} \partial_{x_i} u[\theta] n_j(x)   \left( \int_{\R^d} p_i(v) v_j \, \gamma f \, \dv \right) \d\sigma_{\! x} \\
&=: A+B.
\end{aligned}
$$
Using the decomposition \eqref{eq:fdecomposition},  we get
$$
\int_{\R^d} p_i(v) v_j f \, \dv
= \delta_{ij} \left(1+\frac{2}{d}\right) \theta + \int_{\R^d} p_i(v) v_j f^\perp \, \dv.
$$
As a consequence, we obtain
$$
\begin{aligned}
A 
&=  \left(1+\frac{2}{d}\right) \la -\Delta_x u[\theta] , \theta \ra_{L^2_x(\Omega)} 
- \la \partial_{x_j} \partial_{x_i} u[\theta] , \int_{\R^d} p_i(v) v_j f^\perp \, \dv \ra_{L^2_x(\Omega)}\\
&= \left(1+\frac{2}{d}\right) \|\theta\|^2_{L^2_x(\Omega)} - \la \partial_{x_j} \partial_{x_i} u[\theta] , \int_{\R^d} p_i(v) v_j f^\perp \, \dv \ra_{L^2_x(\Omega)},
\end{aligned}
$$
since by definition of $u[\theta]$ we have $-\Delta_x u[\theta] = \theta$.
Because of~\eqref{eq:uthetafH2}, we obtain
$$
\begin{aligned}
\left| \la \partial_{x_j} \partial_{x_i} u[\theta] , \int_{\R^d} p_i(v) v_j f^\perp \, \dv \ra_{L^2_x(\Omega)}  \right|
&\lesssim \| \nabla_x^2 u[\theta] \|_{L^2_x(\Omega)} \| f^\perp \|_{\HH} \\
&\lesssim \| \theta \|_{L^2_x(\Omega)} \| f^\perp \|_{\HH}.
\end{aligned}
$$
Thanks to Young's inequality, we thus get
$$
A \ge \frac12\left(1+\frac{2}{d}\right) \| \theta \|_{L^2_x(\Omega)}^2 - C \| f^\perp \|_{\HH}^2.  
$$
We now investigate the boundary term $B$. Thanks to Lemma~\ref{lem:boundary}, we have
$$
\begin{aligned}
B &= \int_{\Sigma} \nabla_x u[\theta] \cdot p(v) (\gamma f) \, n(x) \cdot v \, \dv \, \d\sigma_{\! x} \\
&= \int_{\Sigma_{+}} \nabla_x u[\theta] \cdot p(v) \alpha(x) D^\perp f_{+}  \, n(x) \cdot v \, \dv \, \d\sigma_{\! x} \\
&\quad    
+ \int_{\Sigma_{+}} \nabla_x u[\theta] \cdot [p(v) - p(R_x v)] (1-\alpha(x)) D^\perp f_{+}  \, n(x) \cdot v \, \dv \, \d\sigma_{\! x} \\
&\quad    
+ \int_{\Sigma_{+}} \nabla_x u[\theta] \cdot [p(v) - p(R_x v)] D f_{+}  \, n(x) \cdot v \, \dv \, \d\sigma_{\! x} \\
&=: B_1 + B_2 + B_3.
\end{aligned}
$$
We remark that
$$
p(v) - p(R_x v) = 2 n(x) (n(x) \cdot v) \frac{(|v|^2-d-2)}{\sqrt{2d}}
$$
and thus
$$
\nabla_x u[\theta] \cdot [p(v) - p(R_x v)] 
=  2 \nabla_x u[\theta] \cdot n(x) \, (n(x) \cdot v) \, \frac{(|v|^2-d-2)}{\sqrt{2d}}  .
$$
Thanks to the boundary condition satisfied by $u[\theta]$, in the case $\alpha \equiv 0$, we already obtain that $B = 0$. 
Otherwise, when $\alpha \not \equiv 0$, recalling \eqref{eq:def_D}, we first obtain for the term $B_3$, that
$$
\begin{aligned}
B_3 
&= \frac{2c_\mu}{\sqrt{2d}} \int_{\Sigma_+} \nabla_x u[\theta] \cdot n(x) \mu(v) (|v|^2-d-2)\widetilde f(x) \, (n(x) \cdot v)^2 \, \dv \, \d\sigma_{\!x} \\
&= \frac{2c_\mu}{\sqrt{2d}} \int_{\partial\Omega} \nabla_x u[\theta] \cdot n(x) \widetilde f(x) \left( \int_{\Sigma_+^x }  (|v|^2-d-2) \mu(v)  \, (n(x) \cdot v)^2 \dv \right) \d\sigma_{\! x},
\end{aligned}
$$
and the integral in $v$ vanishes, thus $B_3 = 0$.
For the term $B_1$, the Cauchy-Schwarz inequality and~\eqref{eq:uthetafH2} give
$$
\begin{aligned}
|B_1|
&\lesssim \| \nabla_x u[\theta] \|_{L^2_x(\partial\Omega)}  \| \alpha D^\perp f_+ \|_{\partial \HH_+} \\
&\lesssim \| \nabla_x u[\theta] \|_{H^1(\Omega)}  \| \alpha D^\perp f_+ \|_{\partial \HH_+} \\
&\lesssim \| \theta \|_{L^2_x(\Omega)}  \| \alpha D^\perp f_+ \|_{\partial \HH_+} .
\end{aligned}
$$
For the term $B_2$, the boundary condition satisfied by $u[\theta]$ implies
$$
\nabla_x u[\theta] \cdot [p(v) - p(R_x v)] (1-\alpha(x))
= -  \frac{1-\alpha(x)}{2 - \alpha(x)} \alpha(x) u[\theta] 2 (n(x) \cdot v) \frac{(|v|^2-d-2)}{\sqrt{2d}} ,
$$
hence we obtain
$$
\begin{aligned}
|B_2|
&= 2\left| \int_{\Sigma_{+}}   u[\theta] \frac{(|v|^2-d-2)}{\sqrt{2d}}  \alpha(x) \frac{1-\alpha(x)}{2 - \alpha(x)}  D^\perp f_+ \, (n(x) \cdot v)^2 \, \dv \, \d\sigma_{\! x}   \right| \\
&\lesssim \|  u[\theta] \|_{L^2_x(\partial\Omega)}  \| \alpha D^\perp f_+ \|_{\partial \HH_+} \\
&\lesssim \| \theta \|_{L^2_x(\Omega)}  \| \alpha D^\perp f_+ \|_{\partial \HH_+} .
\end{aligned}
$$
We complete the proof by gathering the previous estimates, using Young's inequality and remarking that $\sqrt{\alpha(2-\alpha)} \ge \alpha$.
\end{proof}

\subsection{Momentum}\label{ssec:momentum}

In this subsection we construct a functional that is devised to control the momentum component of the macroscopic part $\pi f$.
We denote 
$$
m[g] := \int_{\R^d}  v g \, \dv , 
$$
so that $m = m[f]$. We define $U[m]$ as the solution to the elliptic equation \eqref{eq:elliptic-korn} associated to $\Xi = m \in L^2_x(\Omega)$ given by Theorem~\ref{theo:regH2-korn}, whence
\begin{equation}\label{eq:UmfH2}
\| U[m] \|_{H^2_x(\Omega)} \lesssim \| m \|_{L^2_x(\Omega)}.
\end{equation}
It is worth noting that in the specular reflection case ($\alpha \equiv 0$ in \eqref{eq:BdyCond}), the condition \eqref{eq:C2C3pif} holds, 
and therefore the solution~$U[m]$ is indeed well-defined. 

\smallskip

Considering the matrix $q_{ij} = (q_{ij})_{1 \le i , j \le d}$ given by 
$$
q_{ij}(v) = v_i v_j - \delta_{ij},
$$
we define the associated moment functional $M_q[g] = (M_{q_{ij}} [g])_{1 \le i , j \le d}$ as
\begin{equation}\label{eq:def-Mq}
M_{q_{ij}}[g] = \int_{\R^d} (v_i v_j - \delta_{ij}) g \, \dv.
\end{equation}

\begin{lem}\label{lem:mLLf}
There holds 
\begin{equation}\label{mLf}
\begin{aligned}
m [\LLL f] 
&= - \nabla_x \varrho - \nabla_x \cdot M_q[f] \\
\end{aligned}
\end{equation}
and 
\begin{equation}\label{Mqf}
M_q[f] = \sqrt{\frac{2}{d}} \theta I_d + M_q[f^\perp].
\end{equation}
As a consequence of Theorem~\ref{theo:regH2-korn}, the unique variational solution $U[m[\LLL f]]$ to \eqref{eq:elliptic-korn-var} associated to $\Xi = m[\LLL f]$ satisfies
\begin{equation}\label{eq:UmLfH1}
\| U[m[\LLL f]] \|_{H^1_x(\Omega)}
\lesssim  
\| \varrho \|_{L^2_x(\Omega)}
+ \| \theta \|_{L^2_x(\Omega)}
+ \| f^\perp \|_{\HH}
+ \| \sqrt{\alpha(2-\alpha)} D^\perp f_+ \|_{\partial \HH_+}.
\end{equation}
\end{lem}

\begin{proof}[Proof of Lemma~\ref{lem:mLLf}]
Writing $\LLL f = - v \cdot \nabla_x f + \CCC f^\perp$ we already obtain that $m[\LLL f] = m [- v \cdot \nabla_x f]$. We hence compute, for $i \in \{1 , \ldots, d \}$,
\begin{equation*}
\begin{aligned}
m_i [- v \cdot \nabla_x f] 
&= - \partial_{x_{j}} \int_{\R^d}  v_i v_j f \, \dv \\
&= - \partial_{x_{i}} \int_{\R^d} f \, \dv 
- \partial_{x_{j}} \int_{\R^d}  (v_i v_j -\delta_{ij}) f \, \dv \\
&= - \partial_{x_i} \varrho - \partial_{x_{j}} M_{q_{ij}}[f],
\end{aligned}
\end{equation*}
which gives \eqref{mLf}. Thanks to the decomposition \eqref{eq:fdecomposition} 
we also obtain, for $i,j \in \{1 , \ldots, d \}$,
\begin{equation*}
\begin{aligned}
M_{q_{ij}}[f]
&= \varrho \int_{\R^d} (v_i v_j - \delta_{ij}) \mu \, \dv 
+m_k   \int_{\R^d} (v_i v_j - \delta_{ij})  v_k \mu \, \dv \\
&\quad 
+\theta \int_{\R^d} (v_i v_j - \delta_{ij})  \left(\frac{|v|^2 - d}{\sqrt{2d}} \right)  \mu \, \dv 
+ M_{q_{ij}}[f^\perp] , 
\end{aligned}
\end{equation*}
which gives \eqref{Mqf} since $\int_{\R^d} (v_i v_j - \delta_{ij}) \mu \, \dv = \int_{\R^d} (v_i v_j - \delta_{ij})  v_k \mu \, \dv= 0$ and $\int_{\R^d} (v_i v_j - \delta_{ij})  \left(\frac{|v|^2 - d}{\sqrt{2d}} \right)  \mu \, \dv = \sqrt{\frac{2}{d}} \delta_{ij}$.

\smallskip

Now let $U := U[m[\LLL f]]$ be the unique variational solution to \eqref{eq:elliptic-korn-var} associated to $\Xi = m[\LLL f]$ from Theorem~\ref{theo:regH2-korn}. 
From Step~1 of the proof of Theorem~\ref{theo:regH2-korn}, one has 
\begin{equation}\label{eq:UmLf1}
\lambda \| U \|_{H^1_x(\Omega)}^2  
\le \| \nabla^s U \|_{L^2_x(\Omega)}^2 
+ \| \sqrt{\tfrac{\alpha}{2-\alpha}} \, U \|_{L^2_x(\partial\Omega)}^2, 
\end{equation}
for some $\lambda >0$. Moreover from \eqref{eq:elliptic-korn-var}, we obtain
\begin{equation}\label{eq:UmLf2}
\begin{aligned}
& \| \nabla^s U \|_{L^2_x(\Omega)}^2 
+ \| \sqrt{\tfrac{\alpha}{2-\alpha}} \, U \|_{L^2_x(\partial\Omega)}^2  \\
&\qquad 
= -   \int_{\Omega} (\nabla_x \varrho + \nabla_x \cdot M_q[f] ) \cdot U \, \dx\\
&\qquad 
=  \int_{\Omega} \varrho I_d : \nabla U \, \dx
+ \int_{\Omega} M_q[f] : \nabla U \, \dx \\
&\qquad \quad
- \int_{\partial \Omega}  \varrho n(x) \cdot U \, \d\sigma_{\! x} 
- \int_{\partial \Omega} M_q[f] n(x) \cdot U  \, \d\sigma_{\! x} \\
&\qquad 
=  \int_{\Omega} \varrho I_d : \nabla^s U \, \dx
+ \int_{\Omega} \left(  \sqrt{\frac{2}{d}} \theta I_d + M_q[f^\perp] \right) : \nabla^s U \, \dx \\
&\qquad \quad
- \int_{\partial \Omega}  M_q[f] n(x) \cdot U \, \d\sigma_{\! x} , 
\end{aligned}
\end{equation}
where we have performed an integration by parts in the second equality, used that $U \cdot n(x) = 0$ since $U \in \VV_\alpha$ and \eqref{Mqf} in the last one. 
We now deal with the boundary term in the last equation. We have, for any $x \in \partial\Omega$,
$$
\begin{aligned}
M_q[f] n(x) \cdot U 
&= \int_{\R^d}  v_i v_j f n_j(x) U_i \, \dv   
- \int_{\R^d}  f n_i(x)  U_i \, \dv  \\
&= \int_{\R^d}  f  (v \cdot U)  (n(x) \cdot v) \, \dv \\
&= \alpha(x) \int_{\Sigma_+^x}    D^\perp f_+ (v \cdot U)   (n(x) \cdot v) \, \dv  \\
&\quad 
+ \int_{\Sigma_+^x}  ( v - R_x v ) \cdot U (1-\alpha(x))  D^\perp f_+  (n(x) \cdot v) \, \dv  \\
&\quad 
+\int_{\Sigma_+^x}   ( v - R_x v ) \cdot U  Df_+  (n(x) \cdot v) \, \dv  , 
\end{aligned}
$$
using that $U \cdot n(x) = 0$ and Lemma~\ref{lem:boundary} in the last line. Observe now that, for any $x \in \partial \Omega$, we have
$$
( v - R_x v ) \cdot U
= 2  \left(n(x) \cdot U \right) (n(x)\cdot v) = 0 ,
$$
by using again that the solution verifies $U \cdot n(x) = 0$. We hence finally get
$$
\begin{aligned}
\left| \int_{\partial\Omega} M_q[f] n(x) \cdot U \, \d\sigma_{\! x}  \right|
&\lesssim \| \sqrt{\alpha(2-\alpha)} \,  D^\perp f_+ \|_{\partial \HH_+} \| \sqrt{\tfrac{\alpha}{2-\alpha}}  U[m[\LLL f]] \|_{L^2_x(\partial\Omega)}.
\end{aligned}
$$

We conclude to \eqref{eq:UmLfH1} by 
gathering this last estimate together with \eqref{eq:UmLf1} and \eqref{eq:UmLf2}, applying Cauchy-Schwarz inequality and remarking that
$$
\| M_q[f^\perp]  \|_{L^2_x(\Omega)}  \lesssim  \| f^\perp  \|_{\HH} .
$$
\end{proof}

We now deduce the following result, which gives a control of the momentum $m$.

\begin{lem}\label{lem:momentum}
There are constants $\kappa_2,C >0$ such that
$$
\begin{aligned}
&\la -\nabla_x^s U[m],  M_q [\LLL f]  \ra_{L^2_x(\Omega)}
+ \la -\nabla_x^s U[m[\LLL f]] ,  M_q [f]  \ra_{L^2_x(\Omega)}\\
&\qquad\qquad
\ge \kappa_2 \| m \|_{L^2_x(\Omega)}^2 
- C \| f^\perp \|_{\HH} \| \varrho \|_{L^2_x(\Omega)}
- C \| \theta \|_{L^2_x(\Omega)} \| \varrho \|_{L^2_x(\Omega)}  \\
&\qquad\qquad\quad
-C \| \theta \|_{L^2_x(\Omega)}^2- C \| f^\perp \|_{\HH}^2 
-C \| \sqrt{\alpha(2-\alpha)} D^\perp f_+ \|_{\partial \HH_+}^2. 
\end{aligned}
$$
\end{lem}

\begin{proof}[Proof of Lemma~\ref{lem:momentum}]
Thanks to \eqref{Mqf} and \eqref{eq:UmLfH1}, we have
$$
\begin{aligned}
&\left|  \la -\nabla_x^s U[m[\LLL f]] ,  \sqrt{\frac{2}{d}} \theta I_d + M_q[f^\perp]  \ra_{L^2_x(\Omega)} \right| \\
&\lesssim \| \nabla_x^s U[m[\LLL f]] \|_{L^2_x(\Omega)} 
\left( \| \theta \|_{L^2_x(\Omega)} + \| f^\perp \|_{\HH} \right)\\
&\lesssim \left(  \| \varrho \|_{L^2_x(\Omega)} + \| \theta \|_{L^2_x(\Omega)} + \| f^\perp \|_{\HH} + \| \sqrt{\alpha(2-\alpha)} D^\perp f_+ \|_{\partial \HH_+}\right) \left( \| \theta \|_{L^2_x(\Omega)} + \| f^\perp \|_{\HH} \right) ,
\end{aligned}
$$
which allows us to bound the second term in the LHS of the estimate of the statement.
For the first term, we write $M_q[\LLL f] = M_q[ - v \cdot \nabla_x f] + M_q[\CCC f^\perp]$ to obtain
$$
\la -\nabla_x^s U[m],  M_q [\LLL f]  \ra_{L^2_x(\Omega)}
= T_1 + T_2,
$$
with 
$$
T_1 := \la (\nabla_x^s U[m])_{ij} , \partial_{x_k} \int_{\R^d} q_{ij}(v) v_k f \, \dv \ra_{L^2_x(\Omega)}
$$
and
$$
T_2 := \la -\nabla_x^s U[m] , \int_{\R^d} q(v) \CCC f^\perp \, \dv \ra_{L^2_x(\Omega)}.
$$
Observing that 
$$
\int_{\R^d} q(v) \CCC f^\perp \, \dv  = \left( f^\perp , \CCC (q \mu) \right)_{L^2_v(\mu^{-1})},
$$
we get from~\eqref{eq:UmfH2} that
$$
|T_2| \lesssim \| \nabla_x^s U[m] \|_{L^2_x(\Omega)} \| f^\perp \|_{\HH}
\lesssim \| m \|_{L^2_x(\Omega)} \| f^\perp \|_{\HH}.
$$
For the term $T_1$, thanks to an integration by parts, we may write 
$$
\begin{aligned}
T_1 
& = - \la \partial_{x_k} (\nabla_x^s U[m])_{ij} ,  \int_{\R^d} q_{ij}(v) v_k f \, \dv \ra_{L^2_x(\Omega)} \\
&\quad 
+ \int_{\partial\Omega}  (\nabla_x^s U[m])_{ij} n_k(x) \left( \int_{\R^d} q_{ij}(v) v_k \, \gamma f \, \dv \right) \d\sigma_{\! x} \\
&=: A+B.  
\end{aligned}
$$
Thanks to the decomposition \eqref{eq:fdecomposition}, we get
$$
\int_{\R^d} q_{ij}(v) v_k f \, \dv
= \delta_{jk} m_i + \delta_{ik} m_j 
+ \int_{\R^d} q_{ij}(v) v_k f^\perp \, \dv, 
$$
and hence 
$$
\begin{aligned}
A 
&=  2 \la - \Div_x ( \nabla_x^s U[m] ) , m \ra_{L^2_x(\Omega)} 
- \la \partial_{x_k} (\nabla_x^s U[m])_{ij} ,  \int_{\R^d} q_{ij}(v) v_k f^\perp \, \dv \ra_{L^2_x(\Omega)} \\
&= 2 \|m\|^2_{L^2_x(\Omega)} - \la \partial_{x_k} (\nabla_x^s U[m])_{ij} ,  \int_{\R^d} q_{ij}(v) v_k f^\perp \, \dv \ra_{L^2_x(\Omega)},
\end{aligned}
$$
since $-\Div_x ( \nabla_x^s U[m] ) = m$ by definition of $U[m]$. Using~\eqref{eq:UmfH2}, we have
$$
\begin{aligned}
\left| \la \partial_{x_k} (\nabla_x^s U[m])_{ij} ,  \int_{\R^d} q_{ij}(v) v_k f^\perp \, \dv \ra_{L^2_x(\Omega)}  \right|
&\lesssim \| \nabla_x^2 U[m] \|_{L^2_x(\Omega)} \| f^\perp \|_{\HH} \\
&\lesssim \| m \|_{L^2_x(\Omega)} \| f^\perp \|_{\HH}.
\end{aligned}
$$
We thus obtain, thanks to Young's inequality,
$$
A \ge \| m \|_{L^2_x(\Omega)}^2 - C \| f^\perp \|_{\HH}^2.  
$$
We now investigate the boundary term $B$. Thanks to Lemma~\ref{lem:boundary}, we have
$$
\begin{aligned}
B &= \int_{\Sigma}  \nabla_x^s U[m] : q(v)  \, \gamma f \, n(x) \cdot v \, \dv \, \d\sigma_{\! x} \\
&= \int_{\Sigma_+}  \nabla_x^s U[m] : q(v)  \alpha(x) D^\perp f_{+} \, n(x) \cdot v \, \dv \, \d\sigma_{\! x} \\
&\quad    
+ \int_{\Sigma_+}  \nabla_x^s U[m] : [q(v)-q(R_x v)] (1-\alpha(x)) D^\perp f_{+} \, n(x) \cdot v \, \dv \, \d\sigma_{\! x} \\
&\quad    
+ \int_{\Sigma_+}  \nabla_x^s U[m] : [q(v)-q(R_x v)] D f_{+} \, n(x) \cdot v \, \dv \, \d\sigma_{\! x} \\
&=: B_1 + B_2 + B_3,
\end{aligned}
$$
and we remark that 
$$
q(v) -  q(R_x v) = 4 \left[ (n(x) \otimes v)^{\operatorname{sym}} -  n(x) \otimes n(x) (n(x) \cdot v)  \right] (n(x) \cdot v),
$$
{ where, for any matrix $M \in \mathcal{M}_d(\R)$, we set $(M^{\operatorname{sym}})_{ij} = \frac12(M_{ij} + M_{ji})$}, so that 
$$
\begin{aligned}
&\nabla_x^s U[m] : [q(v)-q(R_x v)]  \\
&\quad 
= 4 \Big\{ \nabla_x^s U[m] :(n(x) \otimes v)^{\operatorname{sym}}  -  \nabla_x^s U[m] : n(x) \otimes n(x) (n(x) \cdot v) \Big\} (n(x) \cdot v).
\end{aligned}
$$
Taking the scalar product with $v$ in the boundary condition satisfied by $U[m]$, we see that, we already have $B = 0$ in the case $\alpha \equiv 0$.
Otherwise, when $\alpha \not\equiv 0$, we first obtain for the term $B_3$, making a change of variables $v \mapsto R_x v$, using also that $(R_x v \cdot n) = - (v \cdot n)$, and recalling that $D f(x,v) = c_\mu \mu(v) \widetilde f(x)$, that
$$
\begin{aligned}
B_3 
&= 2c_\mu \int_{\Sigma} \nabla_x^s U[m] :  q(v) \mu(v) \widetilde f(x) \, n(x) \cdot v \, \dv \, \d\sigma_{\!x} \\
&= 2c_\mu \int_{\partial\Omega} (\nabla_x^s U[m])_{ij} n_k(x) \widetilde f(x) \left( \int_{\R^d}  q_{ij}(v) v_k \mu(v)  \, \dv \right) \d\sigma_{\! x} = 0,
\end{aligned}
$$
since the integral in $v$ vanishes. 
For the term $B_1$, the Cauchy-Schwarz inequality and~\eqref{eq:UmfH2} give
$$
\begin{aligned}
|B_1|
&\lesssim \| \nabla_x^s U[m] \|_{L^2_x(\partial\Omega)}  \| \alpha D^\perp f_+ \|_{\partial \HH_+} \\
&\lesssim \| m \|_{L^2_x(\Omega)}  \| \alpha D^\perp f_+ \|_{\partial \HH_+} .
\end{aligned}
$$
For the term $B_2$, the boundary condition satisfied by $U[m]$ implies
$$
\begin{aligned}
\nabla_x^s U[m] : [q(v)-q(R_x v)] (1-\alpha(x)) 
&=  - \frac{1-\alpha(x)}{2 - \alpha(x)}  4 \alpha(x) (U[m] \cdot v) (n(x) \cdot v),
\end{aligned}
$$
hence we obtain
$$
\begin{aligned}
|B_2|
&= 4\left|  \int_{\Sigma_{+}}  (U[m] \cdot v) \,  \frac{1-\alpha(x)}{2 - \alpha(x)} \alpha(x) D^\perp f \, (n(x) \cdot v)^2 \, \dv \, \d\sigma_{\! x}   \right|  \\
&\lesssim \|  U[m] \|_{L^2_x(\partial\Omega)}  \| \alpha D^\perp f_+ \|_{\partial \HH_+} \\
&\lesssim \| m \|_{L^2_x(\Omega)}  \| \alpha D^\perp f_+ \|_{\partial \HH_+} .
\end{aligned}
$$
The proof is then complete by gathering previous estimates, using Young's inequality and observing that $\sqrt{\alpha(2-\alpha)} \ge \alpha$ .
\end{proof}

\subsection{Mass}\label{ssec:mass}

In this subsection we introduce the last functional, which is built in order to control the mass component of the macroscopic part $\pi f$.
We denote 
$$
\varrho[g] := \int_{\R^d} g \, \dv , 
$$
so that $\varrho = \varrho[f]$. We consider $u_{\mathrm{N}}[\varrho]$ the solution to the Poisson equation \eqref{eq:elliptic} with Neumann boundary condition associated to $\xi = \varrho \in L^2_x(\Omega)$ constructed in Theorem~\ref{theo:Poisson}, namely $u_{\mathrm{N}} [\varrho]$ satisfies a.e.
\beqn \label{eq:systemuNrho}
\left\{
\begin{aligned}
- \Delta_x u_{\mathrm{N}}[\varrho] &= \varrho \quad \text{in} \quad \Omega , \\
\nabla_x u_{\mathrm{N}}[\varrho] \cdot n(x) &= 0 \quad \text{on} \quad \partial\Omega, 
\end{aligned}
\right.
\eeqn
which is indeed well-defined since $\la \varrho \ra = 0$. In particular, we have
\beqn \label{eq:uNrhoH2}
\| u_{\mathrm{N}}[\varrho] \|_{H^2_x(\Omega)} \lesssim \| \varrho \|_{L^2_x(\Omega)}.
\eeqn

\begin{lem}\label{lem:rLLf}
There holds 
\begin{equation}\label{rhoLf}
\varrho [\LLL f] = -  \nabla_x \cdot m .
\end{equation}
As a consequence of Theorem~\ref{theo:Poisson}, the unique variational solution $u_{\mathrm{N}}[\varrho[\LLL f]]$ to \eqref{eq:varValpha} with Neumann boundary condition associated to $\xi = \varrho[\LLL f]$ satisfies
\begin{equation}\label{eq:uNrhoLfH1}
\| u_{\mathrm{N}}[\varrho[\LLL f]] \|_{H^1_x(\Omega)} \lesssim \| m \|_{L^2_x(\Omega)}.
\end{equation}
\end{lem}

\begin{proof}[Proof of Lemma~\ref{lem:rLLf}]
Since $\LLL f = - v\cdot \nabla_x f + \CCC f^\perp$, one has 
$$
\begin{aligned}
\varrho[\LLL f] 
= \varrho[- v \cdot \nabla_x f] 
= - \nabla_x \cdot \int_{\R^d} v f \, \dv 
\end{aligned}
$$
which gives \eqref{rhoLf}.
Now let $u := u_{\mathrm{N}}[\varrho[\LLL f]]$ be the unique variational solution to \eqref{eq:elliptic} with Neumann boundary condition associated to $\xi = \varrho[\LLL f]$ given by Theorem~\ref{theo:Poisson}. From the variational formulation \eqref{eq:varValpha} we have, thanks to an integration by parts,
$$
\begin{aligned}
\| \nabla_x u \|_{L^2_x(\Omega)}^2 
& = - \int_{\Omega} (\nabla_x \cdot m) u \, \dx \\
&= \int_{\Omega}  m \cdot \nabla_x u \, \dx 
- \int_{\partial\Omega}  m \cdot n(x) \, u \, \d\sigma_{\! x}
= \int_{\Omega}  m \cdot \nabla_x u \, \dx  
\end{aligned}
$$
where we have used that $m \cdot n(x) = 0$ in last equality. We therefore obtain \eqref{eq:uNrhoLfH1} thanks to the Cauchy-Schwarz inequality.
\end{proof}

We now establish the following result, which gives a control of the mass $\varrho$.

\begin{lem}\label{lem:mass}
There are constants $\kappa_3 , C >0$ such that 
$$
\begin{aligned}
&\la -\nabla_x u_{\mathrm{N}}[\varrho],  m [\LLL f]  \ra_{L^2_x(\Omega)}
+ \la -\nabla_x u_{\mathrm{N}}[\varrho[\LLL f]] ,  m[f]  \ra_{L^2_x(\Omega)}\\
&\qquad\qquad
\ge \kappa_3 \| \varrho \|_{L^2_x(\Omega)}^2 
- C \left(\| m \|_{L^2_x(\Omega)}^2 + \| \theta \|_{L^2_x(\Omega)}^2
+ \| f^\perp \|_{\HH}^2 \right) \\ 
&\qquad\qquad\quad 
- C \| { \sqrt{\alpha(2-\alpha)}} D^\perp f_+ \|_{\partial \HH_+}^2 . 
\end{aligned}
$$
\end{lem}

\begin{proof}[Proof of Lemma~\ref{lem:mass}]
From~\eqref{eq:uNrhoLfH1}, we have 
$$
\begin{aligned}
\left| \la -\nabla_x u_{\mathrm{N}}[\varrho[\LLL f]] ,  m[f]  \ra_{L^2_x(\Omega)}   \right| 
&\lesssim \|  \nabla_x u_{\mathrm{N}}[\varrho[\LLL f]] \|_{L^2_x(\Omega)}   \| m \|_{L^2_x(\Omega)} 
\lesssim  \| m \|_{L^2_x(\Omega)}^2,
\end{aligned}
$$
which allows us to bound the second term in the LHS of the estimate of the statement. 
For the first term, writing $m[\LLL f] = m[-v\cdot \nabla_x f] + m[\CCC f^\perp] $ and observing that $m[\CCC f^\perp] = 0$, we obtain
$$
\la -\nabla_x u_{\mathrm{N}}[\varrho],  m [\LLL f]  \ra_{L^2_x(\Omega)} 
= \la \partial_{x_i} u_{\mathrm{N}}[\varrho] , \partial_{x_j} \int_{\R^d} v_i v_j f \, \dv  \ra_{L^2_x(\Omega)}. 
$$
We then write 
$$
\begin{aligned}
&\la \partial_{x_i} u_{\mathrm{N}}[\varrho] , \partial_{x_j} \int_{\R^d} v_i v_j f \, \dv  \ra_{L^2_x(\Omega)} 
\\
&\qquad 
= -\la \partial_{x_j} \partial_{x_i} u_{\mathrm{N}}[\varrho] ,  
\int_{\R^d} v_i v_j f \, \dv  \ra_{L^2_x(\Omega)}
+ \int_{\partial\Omega} \partial_{x_i} u_{\mathrm{N}}[\varrho] n_j(x) 
\left( \int_{\R^d} v_i v_j \, \gamma f \, \dv \right)  \d\sigma_{\! x} \\
&\qquad  =: A+B.
\end{aligned}
$$
Thanks to the decomposition \eqref{eq:fdecomposition}, 
we get
$$
\int_{\R^d} v_i v_j f \, \dv
= \delta_{ij} \varrho + \delta_{ij} \sqrt{\frac{2}{d}} \theta + \int_{\R^d} v_i v_j f^\perp \, \dv,
$$
and hence 
$$
\begin{aligned}
A 
&= \la -\Delta_x u_N[\varrho] , \varrho \ra_{L^2_x(\Omega)} 
+ \sqrt{\frac{2}{d}} \la -\Delta_x u_N[\varrho] , \theta \ra_{L^2_x(\Omega)} 
- \la \partial_{x_j} \partial_{x_i}  u_N[\varrho] , \int_{\R^d} v_i v_j f^\perp \, \dv \ra_{L^2_x(\Omega)}\\
&= \|\varrho\|^2_{L^2_x(\Omega)} + \sqrt{\frac{2}{d}} \la -\Delta_x u_N[\varrho] , \theta \ra_{L^2_x(\Omega)} 
- \la \partial_{x_j} \partial_{x_i}  u_N[\varrho] , \int_{\R^d} v_i v_j f^\perp \, \dv \ra_{L^2_x(\Omega)},
\end{aligned}
$$
since $-\Delta_x u_N[\varrho] = \varrho$ by definition of $u_N[\varrho]$. Using~\eqref{eq:uNrhoH2}, we have
$$
\begin{aligned}
\left| \la \partial_{x_j} \partial_{x_i}  u_N[\varrho] , \int_{\R^d} v_i v_j f^\perp \, \dv \ra_{L^2_x(\Omega)}  \right|
&\lesssim \| \nabla_x^2 u_N[\varrho] \|_{L^2_x(\Omega)} \| f^\perp \|_{\HH} \\
&\lesssim \| \varrho \|_{L^2_x(\Omega)} \| f^\perp \|_{\HH} ,
\end{aligned}
$$
from which it follows, thanks to Young's inequality,
$$
A \ge \frac12 \| \varrho \|_{L^2_x(\Omega)}^2 
- C \| \theta \|_{L^2_x(\Omega)}^2 
- C \| f^\perp \|_{\HH}^2.  
$$
We now investigate the boundary term $B$. Thanks to Lemma~\ref{lem:boundary} we have
$$
\begin{aligned}
B &= \int_{\Sigma} \nabla_x u_{\mathrm{N}}[\varrho] \cdot v \, \gamma f \, n(x) \cdot v \, \dv \, \d\sigma_{\! x} \\
&= \int_{\Sigma_{+}} \nabla_x u_{\mathrm{N}}[\varrho]\cdot v \alpha(x) D^\perp f_{+}  \, n(x) \cdot v \, \dv \, \d\sigma_{\! x} \\
&\quad    
+ \int_{\Sigma_{+}} \nabla_x u_{\mathrm{N}}[\varrho] \cdot [v - R_x v] (1-\alpha(x)) D^\perp f_{+}  \, n(x) \cdot v \, \dv \, \d\sigma_{\! x} \\
&\quad    
+ \int_{\Sigma_{+}} \nabla_x u_{\mathrm{N}}[\varrho] \cdot [v - R_x v] D f_{+}  \, n(x) \cdot v \, \dv \, \d\sigma_{\! x} \\
&=: B_1 + B_2 + B_3,
\end{aligned}
$$
and we remark that
$$
v - R_x v = 2 n(x) (n(x) \cdot v),
$$
so that 
$$
\nabla_x u_{\mathrm{N}}[\varrho] \cdot [v - R_x v] = 2 \nabla_x u_{\mathrm{N}}[\varrho] \cdot n(x) \, (n(x) \cdot v).
$$
Therefore, thanks to the boundary condition satisfied by $u_{\mathrm{N}}[\varrho]$ in~\eqref{eq:systemuNrho}, we already obtain~$B_2 = B_3 = 0$.

In the case $\alpha \equiv 0$, we also have $B_1=0$. Otherwise, when $\alpha \not \equiv 0$, the Cauchy-Schwarz inequality and~\eqref{eq:uNrhoH2} yield 
$$
\begin{aligned}
|B_1|
&\lesssim \| \nabla_x u_{\mathrm{N}}[\varrho] \|_{L^2_x(\partial\Omega)}  \| \alpha D^\perp f_+ \|_{\partial \HH_+} \\
&\lesssim \| \varrho \|_{L^2_x(\Omega)}  \| \alpha D^\perp f_+ \|_{\partial \HH_+} .
\end{aligned}
$$
The proof is then complete by gathering all the previous estimates, using Young's inequality { and observing again that $\sqrt{\alpha(2-\alpha)} \ge \alpha$}.
\end{proof}

\subsection{Proof of Theorem \ref{theo:hypo}}\label{ssec:conclusion}

We define the scalar product $\la\!\la \cdot, \cdot \ra\!\ra$ on~$\HH$ by 
$$
\begin{aligned}
\la \! \la f , g \ra \! \ra 
&:= \la f , g \ra_{\HH} \\
&\quad + \eta_1 \la -\nabla_x u[\theta[f]] ,  M_p [g]  \ra_{L^2_x(\Omega)}
+ \eta_1\la -\nabla_x u[\theta [g]],  M_p [f]  \ra_{L^2_x(\Omega)} \\
&\quad
+ \eta_2 \la -\nabla_x^s U[m[f]] ,  M_q [g]  \ra_{L^2_x(\Omega)}
+ \eta_2 \la -\nabla_x^s U[m[g]] ,  M_q [f]  \ra_{L^2_x(\Omega)} \\
&\quad
+ \eta_3 \la -\nabla_x u_{\mathrm{N}}[\varrho[f]] ,  m [g]  \ra_{L^2_x(\Omega)}
+ \eta_3 \la -\nabla_x u_{\mathrm{N}}[\varrho[g]] ,  m[f]  \ra_{L^2_x(\Omega)}
\end{aligned}
$$
with $0 \ll \eta_3 \ll \eta_2 \ll \eta_1 \ll 1$, and where we recall that the moments $M_p$ and $M_q$ are defined respectively in \eqref{eq:def-Mp} and \eqref{eq:def-Mq}; $u[\theta[f]]$ is the solution of the Poisson equation~\eqref{eq:elliptic} with data $\theta[f]$; $U[m[f]]$ is the solution to the elliptic system \eqref{eq:elliptic-korn} with data $m[f]$; $u_{\mathrm{N}}[\varrho[f]]$ is the solution to the Poisson equation with homogeneous Neumann boundary condition~\eqref{eq:systemuNrho} with data $\varrho[f]$, and similarly for the terms depending on $g$. 
We denote by $\Nt \cdot \Nt$ the norm associated to the scalar product $\la \! \la \cdot , \cdot \ra \! \ra$, and we observe that 
$$
\| f \|_{\HH} \lesssim \Nt f \Nt \lesssim \| f \|_{\HH}.
$$

Let $f$ satisfy the assumptions of Theorem~\ref{theo:hypo}. Recalling that we denote $\varrho=\varrho[f]$, $m=m[f]$ and $\theta=\theta[f]$, noting that $\sqrt{\alpha(2-\alpha)} \ge \alpha$ since $\alpha$ takes values in $[0,1]$, and gathering Lemmas~\ref{lem:micro},~\ref{lem:energy},~\ref{lem:momentum} and~\ref{lem:mass}, one has 
$$
\begin{aligned}
\la \! \la - \LLL f , f \ra \! \ra 
&\ge \lambda \| f^\perp \|_{\HH}^2 + \frac12 \| \sqrt{\alpha(2-\alpha)} D^\perp f_{+} \|_{\partial \HH_+}^2 \\
&\quad 
+\eta_1 \Big( \kappa_1 \| \theta \|_{L^2_x(\Omega)}^2 
- C \| m \|_{L^2_x(\Omega)} \| f^\perp \|_{\HH} 
\\&\qquad \qquad - C \| f^\perp \|_{\HH} ^2 
- C \| \sqrt{\alpha(2-\alpha)} D^\perp f_{+} \|_{\partial \HH_+}^2 \Big)\\
&\quad 
+\eta_2 \Big( \kappa_2 \| m \|_{L^2_x(\Omega)}^2 
- C \| \varrho \|_{L^2_x(\Omega)} \| f^\perp \|_{\HH} 
- C \| \varrho \|_{L^2_x(\Omega)} \| \theta \|_{L^2_x(\Omega)}\\
&\qquad \qquad- C \| \theta \|_{L^2_x(\Omega)}^2 
- C \| f^\perp \|_{\HH}^2 
- C \| \sqrt{\alpha(2-\alpha)} D^\perp f_{+} \|_{\partial \HH_+}^2 \Big) \\
&\quad 
+\eta_3 \Big( \kappa_3 \| \varrho \|_{L^2_x(\Omega)}^2 
- C \| m \|_{L^2_x(\Omega)}^2 
- C \| \theta \|_{L^2_x(\Omega)}^2 \\
&\qquad \qquad- C \| f^\perp \|_{\HH}^2 
- C \| \sqrt{\alpha(2-\alpha)} D^\perp f_{+} \|_{\partial \HH_+}^2  \Big).
\end{aligned}
$$
Thanks to Young's inequality,  we have
$$
\begin{aligned}
\eta_1 C \| m \|_{L^2_x(\Omega)} \| f^\perp \|_{\HH}
&\le  \frac{\lambda}{4} \| f^\perp \|_{\HH}^2   + C \eta_1^2 \| m \|_{L^2_x(\Omega)}^2, \\
\eta_2 C \| \varrho \|_{L^2_x(\Omega)} \| f^\perp \|_{\HH} 
&\le \frac{\lambda}{4} \| f^\perp \|_{\HH}^2 + C \eta_2^2 \| \varrho \|_{L^2_x(\Omega)}^2 ,
\\
\eta_2 C \| \varrho \|_{L^2_x(\Omega)} \| \theta \|_{L^2_x(\Omega)} 
&\le \frac{\eta_1 \kappa_1}{2} \| \theta \|_{L^2_x(\Omega)}^2 + C \frac{\eta_2^2}{\eta_1} \| \varrho \|_{L^2_x(\Omega)}^2  .
\end{aligned}
$$
We thus obtain
$$
\begin{aligned}
\la \! \la - \LLL f , f \ra \! \ra 
&\ge \left(\frac{\lambda}{2} - \eta_1 C - \eta_2 C - \eta_3 C\right) \| f^\perp \|_{\HH}^2 \\
&\quad 
+ \left(\frac12 - \eta_1 C - \eta_2 C - \eta_3 C  \right) \| \sqrt{\alpha(2-\alpha)} D^\perp f_{+} \|_{\partial \HH_+}^2 \\
&\quad 
+\left(\frac{\eta_1 \kappa_1}{2} - \eta_2 C - \eta_3 C \right)\| \theta \|_{L^2_x(\Omega)}^2  \\
&\quad 
+ \left( \eta_2 \kappa_2 - \eta_1^2 C - \eta_3 C \right)\| m \|_{L^2_x(\Omega)}^2  \\
&\quad 
+ \left( \eta_3 \kappa_3  - \eta_2^2 C - \frac{\eta_2^2}{\eta_1}C\right)\| \varrho \|_{L^2_x(\Omega)}^2 .
\end{aligned}
$$
We now choose $\eta_1 := \eta$, $\eta_2 := \eta^{\frac{3}{2}}$, $\eta_3 := \eta^{\frac{7}{4}}$, and we deduce
$$
\begin{aligned}
\la \! \la - \LLL f , f \ra \! \ra 
&\ge \left(\frac{\lambda}{2} - \eta C \right) \| f^\perp \|_{\HH}^2 
+ \left(\frac12 - \eta C   \right) \| \sqrt{\alpha(2-\alpha)} D^\perp f_{+} \|_{\partial \HH_+}^2 \\
&\quad 
+\eta \left( \frac{\kappa_1}{2} - \eta^{\frac{1}{2}} C \right)\| \theta \|_{L^2_x(\Omega)}^2   
+ \eta^{\frac{3}{2}} \left(  \kappa_2 - \eta^{\frac{1}{4}} C  \right)\| m \|_{L^2_x(\Omega)}^2  \\
&\quad 
+ \eta^{\frac{7}{4}} \left(  \kappa_3  - \eta^{\frac{1}{4}} C\right)\| \varrho \|_{L^2_x(\Omega)}^2 .
\end{aligned}
$$
Choosing $0 < \eta < 1$ small enough, we get
$$
\begin{aligned}
\la \! \la - \LLL f , f \ra \! \ra 
&\ge \kappa \left( \| f^\perp \|_{\HH}^2 
+\| \varrho \|_{L^2_x(\Omega)}^2 
+\| m \|_{L^2_x(\Omega)}^2 
+\| \theta \|_{L^2_x(\Omega)}^2  \right) \\
&\quad 
+ \kappa' \| \sqrt{\alpha(2-\alpha)} D^\perp f_{+} \|_{\partial \HH_+}^2
\end{aligned}
$$
for some constants $\kappa,\kappa' >0$.
We conclude the proof of Theorem~\ref{theo:hypo} since 
$$
\| f^\perp \|_{\HH}^2 
+\| \varrho \|_{L^2_x(\Omega)}^2 
+\| m \|_{L^2_x(\Omega)}^2 
+\| \theta \|_{L^2_x(\Omega)}^2 = \| f \|_{\HH}^2
$$
and $\| \cdot \|_{\HH}$ is equivalent to $\Nt \cdot \Nt$. 
\qed

\section{Weakly coercive operators}
\label{sec:weak}

In this section we extend our method to the case in which the collision operator $\CCC$ is \emph{weakly coercive}, that is, it satisfies assumption (A2') below which is weaker than the coercive estimate of assumption (A2) in Subsection~\ref{subsec:equation}. 

In this situation we do not expect to obtain an exponential decay but only a sub-exponential decay supposing further integrability/regularity properties of the initial data; in other words the semigroup associated to the full linear operator $\LLL$ is not uniformly exponentially stable but only strongly stable.

These weakly coercive operators arise naturally in several classes of evolution PDEs. 
In the setting of control theory and wave-type equations we refer to the works \cite{Lebeau,LebeauRobb,Burq,MR3219503,MR3626005} and the references therein, in which the energy of the equation is shown to decay with non-exponential rate. These results have then inspired an abstract theory for strongly stable semigroups. We refer to \cite{BD,BEPS,BCT} and the references therein, where such a line of research is developped. 

In the framework of kinetic equations, the works \cite{Caflisch1,Caflisch2} have established the sub-exponential decay of the semigroup associated to the linearized cutoff Boltzmann equation with soft potentials. We also refer to the works \cite{GS1,GS2} that establish decay estimates for the non-cutoff Boltzmann and Landau equations with very soft potentials, as well as \cite{MR3625186} for the Landau equation. All these results are established in the torus or the whole space, and, to the best of our knowledge, the only works concerning domains with boundary conditions are the recent results of \cite{MR4076068} for the Landau equation with specular reflection boundary condition, and \cite{DuanLiuSakamotoStrain} for non-cutoff Boltzmann and Landau equations in a finite channel with specular reflection or inflow boundary conditions.
Concerning Fokker-Planck equations  and kinetic Fokker-Planck equations we shall quote \cite{RockWang,KM} and \cite{MR4069622}, as well as the references therein. We also mention the results concerning degenerate linear transport equations \cite{MR2569870,MR3048598,MR4063917}, as well as degenerate linear Boltzmann equations \cite{MR3479064}. Finally, the free transport equation with diffusive or Maxwell boundary condition has been tackled in \cite{MR2765738,MR3199988,MR4179249} for instance. 

\smallskip
We assume in this section that the operator $\CCC$ satisfies (A1) on $L^2_v (\mu^{-1})$, as well as: 
\begin{itemize}


\item[(A2')] The operator is self-adjoint on $L^2_v(\mu^{-1})$ and negative $(\CCC f , f )_{L^2_v(\mu^{-1})} \le 0$, so that its spectrum is included in $\R_{-}$, and \eqref{eq:local-conservations} holds true for any { $g \in \mathrm{Dom}(\CCC)$}. We assume further that $\CCC$ satisfies a weak coercivity estimate: there is a positive constant $\lambda >0$ and a radially symmetric function $\omega_0 : \R^d \to [1,\infty)$ with $\lim_{|v| \to \infty} \omega_0(v) = \infty$ such that for any $f \in \mathrm{Dom} (\CCC)$ one has
$$
(-\CCC f , f)_{L^2_v(\mu^{-1})} \ge \lambda \|  f^\perp \|_{L^2_{v}( \omega_0^{-1} \mu^{-1})}^2,
$$
where $f^\perp := f - \pi f$. 

\item[(A3')] For any polynomial function $\phi=\phi(v) : \R^d \to \R$ of degree $\le 4$, one has $ \mu \phi \in \mathrm{Dom} (\CCC)$ with
$$
\| \CCC (\phi \mu) \|_{L^2_v(\omega_0 \mu^{-1})}  < \infty,
$$
and, for some positive constant $C > 0$, for all $f \in \mathrm{Dom}(\CCC)$,
$$
\left|  \int_{\R^d} \phi(v) f^\perp \, \dv \right| \le C \| f^\perp \|_{L^2_{v}(\omega_0^{-1}\mu^{-1})}.
$$


\item[(A4)] There exists a radially symmetric function $\omega_1 : \R^d \to [1,\infty)$ with $\lim_{|v| \to \infty} \omega_1(v) = \infty$ and a positive constant $C>0$ such that for any $f \in \mathrm{Dom} (\LLL)$,
 one has
$$
\la  \LLL f , f \ra_{L^2_{x,v} ( \omega_1 \mu^{-1})} \le C \| f \|_{L^2_{x,v}(\omega_0^{-1} \mu^{-1})}^2. 
$$

\end{itemize}

\medskip
We recall that $\HH = L^2_{x,v}(\mu^{-1})$ and in this Section, we will also use the following notations: $\HH_0 := L^2_{x,v} (\omega_0^{-1} \mu^{-1})$ and $\HH_1 := L^2_{x,v} (\omega_1 \mu^{-1})$.
Remark now that we have 
$$
\|f ^{\perp}\|^2_{\HH_0} + \| \pi f \|^2_{\HH_0} \lesssim \| f \|_{\HH_0}^2 \lesssim  
\|f ^{\perp}\|^2_{\HH_0} + \| \pi f \|^2_{\HH_0}
$$
and
$$ 
\|\pi f\|^2_{\HH_0} \lesssim \| \varrho \|^2_{L^2_x(\Omega)} + \| m \|^2_{L^2_x(\Omega)} + \| \theta \|_{L^2_x(\Omega)}^2 \lesssim \|\pi f\|^2_{\HH_0}. 
$$

Repeating the proof of Theorem~\ref{theo:hypo} with the above assumptions we obtain:
\begin{theo}\label{theo:weak-hypo}
There exists a scalar product $\la\!\la \cdot , \cdot \ra \! \ra_{\HH}$ on the space $\HH$ so that the associated norm $\Nt \cdot \Nt_{\HH}$ is equivalent to the usual norm $\| \cdot \|_{\HH}$, and for which the linear operator $\LLL$ satisfies the following weak coercivity estimate: there is a positive constant $\kappa >0$ such that 
one has 
$$
\la \! \la - \LLL f , f \ra\!\ra_{\HH} \ge \kappa  \| f \|_{\HH_0}^2
$$
for any $f \in \mathrm{Dom}(\LLL)$ satisfying the boundary condition \eqref{eq:BdyCond}, assumption~\eqref{eq:C1} and furthermore assumptions \eqref{eq:C2}-\eqref{eq:C3} in the specular reflection case ($\alpha\equiv 0$ in \eqref{eq:BdyCond}).
%
%

\end{theo}

As a consequence of the weak coercivity estimate for $\LLL$, we obtain the following result of sub-exponential decay to equilibrium.

\begin{theo}\label{theo:weak-main}
Let $f_{\mathrm{in}} \in \HH_1$ satisfying 
condition~\eqref{eq:C1} and furthermore~\eqref{eq:C2}-\eqref{eq:C3} in the specular reflection case ($\alpha\equiv 0$ in \eqref{eq:BdyCond}).
There exist a positive constant~$C>0$ and a decreasing function $\vartheta : \R_+ \to \R_+$ with $\lim_{t \to \infty} \vartheta(t) = 0$ such that for any solution $f$ to \eqref{eq:dtf=Lf}--\eqref{eq:BdyCond} (with $\CCC$ satisfying (A1)--(A2')--(A3')--(A4) above) associated to the initial data~$f_{\mathrm{in}}$,  there holds
$$
\| f(t) \|_{\HH} \le C \vartheta(t) \| f_{\mathrm{in}} \|_{\HH_1}, \quad \forall \, t \ge 0.
$$
\end{theo}

\begin{proof}[Proof of Theorem~\ref{theo:weak-main}]
Let $f$ be a solution to \eqref{eq:dtf=Lf}--\eqref{eq:BdyCond} associated to  $f_{\mathrm{in}} \in \mathrm{Dom}(\LLL)$, the general case when  $f_{\mathrm{in}} \in \HH_1$ then deduces by a usual density argument. Thanks to Theorem~\ref{theo:weak-hypo}, we have 
\beqn\label{eq:weak1}
\frac{\d}{\dt} \Nt f(t) \Nt_{\HH}^2 
= \la \! \la \LLL f(t) , f(t) \ra\!\ra_{\HH} 
\le -\kappa  \| f(t) \|_{\HH_0}^2.
\eeqn
Remark that for any $R>0$ we have the following interpolation inequality
\beqn\label{eq:weak2}
\| g \|^2_{\HH}
\le \omega_0(R) \| g \|^2_{\HH_0}
+ \frac{1}{\omega_1(R)} \, \| g \|^2_{\HH_1}.
\eeqn
Moreover we claim that there is a constant $C>0$ such that 
\beqn\label{eq:weak3}
\| f(t) \|_{\HH_1} \le C \| f_{\mathrm{in}} \|_{\HH_1}.
\eeqn
Indeed for $\delta >0$ small enough, we define the following scalar product on $\HH_1$ 
$$
\la \! \la f , g  \ra \! \ra_{\HH_1} 
:= \delta  \la  f , g \ra_{\HH_1} + \la\!\la   f , g \ra\!\ra_{\HH}.
$$
Gathering (A4) and Theorem~\ref{theo:weak-hypo}, we obtain
$$
\begin{aligned}
\la \! \la \LLL f , f \ra  \! \ra_{\HH_1} 
&\le (\delta C - \kappa) \| f \|_{\HH_0}^2 \le 0,
\end{aligned}
$$
which implies the claim by observing that the norm associated to $ \la \! \la \cdot , \cdot \ra \! \ra_{\HH_1} $ is equivalent to the standard norm on $\HH_1$.

From \eqref{eq:weak1}, \eqref{eq:weak2} and \eqref{eq:weak3}, we therefore deduce
$$
\begin{aligned}
\frac{\d}{\dt} \Nt f(t) \Nt_{\HH}^2 
&\le -\frac{\kappa}{\omega_0(R)} \, \| f(t) \|_{\HH}^2
+ \frac{\kappa}{\omega_0(R) \omega_1(R)} \, \| f(t) \|_{\HH_1}^2 \\
&\le -\frac{c \kappa}{\omega_0(R)} \, \Nt f(t) \Nt_{\HH}^2
+ \frac{\kappa C}{\omega_0(R) \omega_1(R)} \,\| f_{\mathrm{in}} \|_{\HH_1}^2,
\end{aligned}
$$
 for some constant $c>0$,
where we have used in last line that $\Nt \cdot \Nt_{\HH}$ and $\| \cdot \|_{\HH}$ are equivalent, and the above claim.
From the above inequality it follows
$$
\begin{aligned}
\Nt f(t) \Nt_{\HH}^2 
&\le \exp\left( -\frac{c \kappa}{\omega_0(R)} \, t \right) \Nt f_{\mathrm{in}} \Nt_{\HH}^2 
+\frac{C}{c \omega_1(R)} \, \| f_{\mathrm{in}} \|_{\HH_1}^2 \\
&\le \left\{ \exp\left( -\frac{c \kappa}{\omega_0(R)} \, t \right) +  \frac{C}{c \omega_1(R)}  \right\} \| f_{\mathrm{in}} \|_{\HH_1}^2,
\end{aligned}
$$
for any $R>0$. 
Defining
$$
\vartheta(t) := \left( \inf_{R>0} \left\{ \exp\left( -\frac{c \kappa}{\omega_0(R)} \, t \right) +  \frac{C}{c \omega_1(R)}  \right\} \right)^{\frac12}, 
$$
we hence obtain 
$$
\Nt f(t) \Nt_{\HH} \le  \vartheta(t) \| f_{\mathrm{in}} \|_{\HH_1},
$$
which concludes the proof using again that $\Nt \cdot \Nt_{\HH}$ and $\| \cdot \|_{\HH}$ are equivalent.
\end{proof}

\section{Hydrodynamic limits}
\label{sec:hydro}

In this part, we study the following rescaled problem: 
\bear\label{eq:dtf=Lepsf}
\partial_t f &=& \LLL_\eps f := - {1 \over \eps} v \cdot \nabla_x f +  {1 \over \eps^2} \CCC f \quad\hbox{in}\quad (0,\infty) \times \OO, 
\\ \label{eq:BdyCondeps}
\gamma_{\!-} f &=& \RRR \gamma_{\!+}  f \quad\hbox{on}\quad (0,\infty) \times \Sigma, 
\eear
with $\eps \in (0,1]$, $\CCC$ satisfying assumptions~(A1),~(A2) and~(A3) introduced in Subsection~\ref{subsec:equation} and the boundary condition~\eqref{eq:BdyCondeps} being the same as~\eqref{eq:BdyCond} described in Subsection~\ref{subsec:equation}. The motivation to study this problem comes from the issue of deriving the incompressible Navier-Stokes-Fourier system from kinetic equations. Indeed, it is well-known (see~\cite{BGL1}) that in order to reach this goal, we shall introduce the dimensionless Knudsen number $\eps$ and the problem reduces to the analysis of the following equation
\beqn \label{eq:nonlineareps}
\partial_t F^\eps = \LLL_\eps F^\eps + {1 \over \eps} Q(F^\eps,F^\eps)
\eeqn
with 
$$
\LLL_\eps f = - {1 \over \eps} v \cdot \nabla_x f +  {1 \over \eps^2} \CCC f, \quad \CCC f := Q(\mu,f) + Q(f,\mu). 
$$
Then, in order to derive the incompressible Navier-Stokes-Fourier limit from kinetic equations, the purpose is to prove that, as $\eps$ goes to $0$, a solution $F^\eps$ to~\eqref{eq:nonlineareps} converges towards some limit that depends on time and space variables only through macroscopic quantities that are solutions to the incompressible Navier-Stokes-Fourier system. The starting point of this study is the analysis of the linearized problem~\eqref{eq:dtf=Lepsf} and our method is robust enough to treat this rescaled problem. More precisely, we are able to provide a result of large time stability for the linear problem~\eqref{eq:dtf=Lepsf}-\eqref{eq:BdyCondeps} uniformly with respect to the parameter~$\eps >0$.

\medskip
The problem of deriving incompressible Navier-Stokes equation from Boltzmann equation has been largely studied in the framework of weak solutions (renormalized for the Boltzmann equation and Leray type for the Navier-Stokes one), in the torus, the whole space or bounded domains. We do not make an extensive presentation here of this type of result but just mention the papers~\cite{BGL1,BGL2} in which this program has been initiated and~\cite{GSR} in which the first complete proof of convergence has been obtained in the whole space. We also mention the works~\cite{MSR,SR-book,JM} in which the problem has been treated in bounded domains starting from  the renormalized solutions constructed in~\cite{MischlerENS2010}.

Concerning the case of strong solutions, we mention the works~\cite{BardosUkai,DMEL} and more recent ones~\cite{ALT,Briant,BMM,GT,JXZ,Rachid} which are all framed in the torus and/or the whole space. To our knowledge, no result of derivation is available for strong solutions in a bounded domain. The study of this derivation will be the object of a forthcoming work. We focus here on the study of the linearized rescaled problem~\eqref{eq:dtf=Lepsf}-\eqref{eq:BdyCondeps}.  

\medskip
We here give an adapted version of Theorem~\ref{theo:hypo} in our new rescaled framework:
\begin{theo}\label{theo:hydro-hypo}
There exists a scalar product $\la\!\la \cdot , \cdot \ra \! \ra_{\!\eps}$ on the space $\HH$ so that the associated norm $\Nt \cdot \Nt_\eps$ is equivalent to the usual norm $\| \cdot \|_{\HH}$ uniformly in $\eps \in (0,1]$, and for which the linear operator $\LLL_\eps$ satisfies the following coercivity estimate: there is a positive constant $\kappa >0$ such that for any $\eps \in (0,1]$,
one has 
$$
\la \! \la - \LLL_\eps f , f \ra\!\ra_{\!\eps} \ge \kappa  \Nt f \Nt^2_\eps + {\kappa \over \eps^2} \|f^\perp\|_{\HH}^2, 
$$
for any $f \in \mathrm{Dom}(\LLL)$ satisfying the boundary condition \eqref{eq:BdyCondeps}, assumption~\eqref{eq:C1} and furthermore assumptions \eqref{eq:C2}-\eqref{eq:C3} in the specular reflection case ($\alpha\equiv 0$ in \eqref{eq:BdyCond}).
\end{theo}

\begin{proof}[Sketch of the proof of Theorem~\ref{theo:hydro-hypo}]
Using the same notations as in Subsection~\ref{ssec:conclusion}, we introduce the following scalar product on $\HH$: 
$$
\begin{aligned}
\la \! \la f , g \ra \! \ra_{\!\eps}
&:= \la f , g \ra_{\HH} \\
&\quad + \eta_1 \eps \la -\nabla_x u[\theta[f]] ,  M_p [g]  \ra_{L^2_x(\Omega)} 
+ \eta_1 \eps\la -\nabla_x u[\theta [g]],  M_p [f]  \ra_{L^2_x(\Omega)} \\
&\quad
+ \eta_2 \eps \la -\nabla_x^s U[m[f]] ,  M_q [g] \ra_{L^2_x(\Omega)}
+\eta_2 \eps\la -\nabla_x^s U[m[g]] ,  M_q [f]  \ra_{L^2_x(\Omega)} \\
&\quad
+ \eta_3 \eps \la -\nabla_x u_{\mathrm{N}}[\varrho[f]] ,  m [g]  \ra_{L^2_x(\Omega)}
+ \eta_3 \eps \la -\nabla_x u_{\mathrm{N}}[\varrho[g]] ,  m[f]  \ra_{L^2_x(\Omega)}
\end{aligned}
$$
with $0 \ll \eta_3 \ll \eta_2 \ll \eta_1 \ll 1$ chosen as in the proof of Theorem~\ref{theo:hypo}.
We denote by $\Nt \cdot \Nt_\eps$ the norm associated to the scalar product~$\la \! \la \cdot , \cdot \ra \! \ra_{\!\eps}$, and we observe that 
$$
\| f \|_{\HH} \lesssim \Nt f \Nt_\eps \lesssim \| f \|_{\HH}
$$
where the multiplicative constants are uniform in $\eps \in (0,1]$. The norms $\| \cdot \|_{\HH}$ and $\Nt \cdot \Nt_\eps$ are thus equivalent independently of $\eps \in (0,1]$. Repeating the proof of Theorem~\ref{theo:hypo}, we obtain the desired result. 
\end{proof}

\medskip

Using once more this equivalence of norms, we are able to prove the following stability result for our equation~\eqref{eq:dtf=Lepsf}-\eqref{eq:BdyCondeps} uniformly in~$\eps \in (0,1]$:
\begin{theo}\label{theo:hydro-main}
Let $f_{\mathrm{in}}^\eps \in \HH$ satisfying condition~\eqref{eq:C1} and furthermore \eqref{eq:C2}-\eqref{eq:C3} in the specular reflection case ($\alpha\equiv 0$ in \eqref{eq:BdyCond}).
There exist positive constants $\kappa , C >0$ independent of $\eps \in (0,1]$ such that for any solution~$f^\eps$ to \eqref{eq:dtf=Lepsf}--\eqref{eq:BdyCondeps} associated to the initial data $f_{\mathrm{in}}^\eps$, for any $\eps \in (0,1]$ and for any  $t \ge 0$, there holds
$$
\| f^\eps(t) \|_{\HH}  \le C e^{-\kappa t} \| f_{\mathrm{in}}^\eps \|_{\HH} \quad \text{and} \quad  {1 \over \eps^2} \int_0^\infty \|(f^\eps)^\perp(s)\|^2_{\HH} \, e^{2\kappa s} \, \d s \le C \| f_{\mathrm{in}}^\eps \|^2_{\HH}.  
$$
\end{theo}

\begin{rem}
Notice that we can perform the same analysis to extend the result of this Section to operators that satisfy a weak coercivity estimate as in Section~\ref{sec:weak}. Namely, one can obtain sub-exponential decay of the solution $f^\eps$ to~\eqref{eq:dtf=Lepsf}-\eqref{eq:BdyCondeps}, that is uniform in $\eps \in (0,1]$, when the collision operator $\CCC$ involved in~\eqref{eq:dtf=Lepsf} satisfy assumptions (A1), (A2'), (A3') and (A4) of Section~\ref{sec:weak}. 
\end{rem}


\phantomsection
\bibliographystyle{acm}
\addcontentsline{toc}{section}{References}

\begin{thebibliography}{10}

\bibitem{ALT}
{\sc Alonso, R.~J., Lods, B., and Tristani, I.}
\newblock Fluid dynamic limit of {B}oltzmann equation for granular hard-spheres
  in a nearly elastic regime.
\newblock Preprint arXiv:2008.05173.

\bibitem{MR3219503}
{\sc Anantharaman, N., and L\'{e}autaud, M.}
\newblock Sharp polynomial decay rates for the damped wave equation on the
  torus.
\newblock {\em Anal. PDE 7}, 1 (2014), 159--214.
\newblock With an appendix by St\'{e}phane Nonnenmacher.

\bibitem{MR2765738}
{\sc Aoki, K., and Golse, F.}
\newblock On the speed of approach to equilibrium for a collisionless gas.
\newblock {\em Kinet. Relat. Models 4}, 1 (2011), 87--107.

\bibitem{ArkerydMaslova}
{\sc Arkeryd, L., and Maslova, N.}
\newblock On diffuse reflection at the boundary for the {B}oltzmann equation
  and related equations.
\newblock {\em J. Statist. Phys. 77}, 5-6 (1994), 1051--1077.

\bibitem{MR274925}
{\sc Bardos, C.}
\newblock Probl\`emes aux limites pour les \'{e}quations aux d\'{e}riv\'{e}es
  partielles du premier ordre \`a coefficients r\'{e}els; th\'{e}or\`emes
  d'approximation; application \`a l'\'{e}quation de transport.
\newblock {\em Ann. Sci. \'{E}cole Norm. Sup. (4) 3\/} (1970), 185--233.

\bibitem{BGL2}
{\sc Bardos, C., Golse, F., and Levermore, C.~D.}
\newblock Fluid dynamic limits of kinetic equations. {II}. {C}onvergence proofs
  for the {B}oltzmann equation.
\newblock {\em Comm. Pure Appl. Math. 46}, 5 (1993), 667--753.

\bibitem{BGL1}
{\sc Bardos, C., Golse, F., and Levermore, D.}
\newblock Fluid dynamic limits of kinetic equations. {I}. {F}ormal derivations.
\newblock {\em J. Statist. Phys. 63}, 1-2 (1991), 323--344.

\bibitem{BardosUkai}
{\sc Bardos, C., and Ukai, S.}
\newblock The classical incompressible {N}avier-{S}tokes limit of the
  {B}oltzmann equation.
\newblock {\em Math. Models Methods Appl. Sci. 1}, 2 (1991), 235--257.

\bibitem{BEPS}
{\sc B{\'a}tkai, A., Engel, K.-J., Pr{\"u}ss, J., and Schnaubelt, R.}
\newblock Polynomial stability of operator semigroups.
\newblock {\em Math. Nachr. 279}, 13-14 (2006), 1425--1440.

\bibitem{BCT}
{\sc Batty, C.~J.~K., Chill, R., and Tomilov, Y.}
\newblock Fine scales of decay of operator semigroups.
\newblock {\em J. Eur. Math. Soc. 18}, 4 (2016), 853--929.

\bibitem{BD}
{\sc Batty, C. J.~K., and Duyckaerts, T.}
\newblock Non-uniform stability for bounded semi-groups on {B}anach spaces.
\newblock {\em J. Evol. Equ. 8}, 4 (2008), 765--780.

\bibitem{MR3048598}
{\sc Bernard, E., and Salvarani, F.}
\newblock On the convergence to equilibrium for degenerate transport problems.
\newblock {\em Arch. Ration. Mech. Anal. 208}, 3 (2013), 977--984.

\bibitem{MR4179249}
{\sc Bernou, A.}
\newblock A semigroup approach to the convergence rate of a collisionless gas.
\newblock {\em Kinet. Relat. Models 13}, 6 (2020), 1071--1106.

\bibitem{MR2150445}
{\sc Boyer, F.}
\newblock Trace theorems and spatial continuity properties for the solutions of
  the transport equation.
\newblock {\em Differential Integral Equations 18}, 8 (2005), 891--934.

\bibitem{Briant}
{\sc Briant, M.}
\newblock From the {B}oltzmann equation to the incompressible {N}avier-{S}tokes
  equations on the torus: a quantitative error estimate.
\newblock {\em J. Differential Equations 259}, 11 (2015), 6072--6141.

\bibitem{MR3579575}
{\sc Briant, M.}
\newblock Perturbative theory for the {B}oltzmann equation in bounded domains
  with different boundary conditions.
\newblock {\em Kinet. Relat. Models 10}, 2 (2017), 329--371.

\bibitem{MR3562318}
{\sc Briant, M., and Guo, Y.}
\newblock Asymptotic stability of the {B}oltzmann equation with {M}axwell
  boundary conditions.
\newblock {\em J. Differential Equations 261}, 12 (2016), 7000--7079.

\bibitem{BMM}
{\sc Briant, M., Merino-Aceituno, S., and Mouhot, C.}
\newblock From {B}oltzmann to incompressible {N}avier-{S}tokes in {S}obolev
  spaces with polynomial weight.
\newblock {\em Anal. Appl. (Singap.) 17}, 1 (2019), 85--116.

\bibitem{Burq}
{\sc Burq, N.}
\newblock D\'ecroissance de l'\'energie locale de l'\'equation des ondes pour
  le probl\`eme ext\'erieur et absence de r\'esonance au voisinage du r\'eel.
\newblock {\em Acta Math. 180}, 1 (1998), 1--29.

\bibitem{MR4063917}
{\sc Ca\~{n}izo, J.~A., Cao, C., Evans, J., and Yolda\c{s}, H.}
\newblock Hypocoercivity of linear kinetic equations via {H}arris's theorem.
\newblock {\em Kinet. Relat. Models 13}, 1 (2020), 97--128.

\bibitem{Caflisch1}
{\sc Caflisch, R.~E.}
\newblock The {B}oltzmann equation with a soft potential. {I}. {L}inear,
  spatially-homogeneous.
\newblock {\em Comm. Math. Phys. 74}, 1 (1980), 71--95.

\bibitem{Caflisch2}
{\sc Caflisch, R.~E.}
\newblock The {B}oltzmann equation with a soft potential. {II}. {N}onlinear,
  spatially-periodic.
\newblock {\em Comm. Math. Phys. 74}, 2 (1980), 97--109.

\bibitem{MR4069622}
{\sc Cao, C.}
\newblock The kinetic {F}okker-{P}lanck equation with weak confinement force.
\newblock {\em Commun. Math. Sci. 17}, 8 (2019), 2281--2308.

\bibitem{CDHMMS}
{\sc Carrapatoso, K., Dolbeault, J., H\'erau, F., Mischler, S., Mouhot, C., and
  Schmeiser, C.}
\newblock Linear stability of a confined system of charged particles.
\newblock In preparation.

\bibitem{MR3625186}
{\sc Carrapatoso, K., and Mischler, S.}
\newblock Landau equation for very soft and {C}oulomb potentials near
  {M}axwellians.
\newblock {\em Ann. PDE 3}, 1 (2017), Paper No. 1, 65.

\bibitem{MR777741}
{\sc Cessenat, M.}
\newblock Th\'{e}or\`emes de trace pour des espaces de fonctions de la
  neutronique.
\newblock {\em C. R. Acad. Sci. Paris S\'{e}r. I Math. 300}, 3 (1985), 89--92.

\bibitem{MR2119999}
{\sc Ciarlet, P.~G., and Ciarlet, Jr., P.}
\newblock Another approach to linearized elasticity and a new proof of {K}orn's
  inequality.
\newblock {\em Math. Models Methods Appl. Sci. 15}, 2 (2005), 259--271.

\bibitem{CostabelDN}
{\sc Costabel, M., Dauge, M., and Nicaise, S.}
\newblock Corner singularities and analytic regularity for linear elliptic
  systems. part i: Smooth domains.
\newblock Pr{\'e}publication IRMAR, 10-09, hal-00453934v2.

\bibitem{DGineq}
{\sc Darroz\`es, J.~S., and Guiraud, J.~P.}
\newblock G{\'e}n{\'e}ralisation formelle du {T}h{\'e}or{\`e}me {H} en
  pr{\'e}sence de parois.
\newblock {\em C. R. Acad. Sci. Paris 262\/} (1966), 368--371.

\bibitem{DMEL}
{\sc De~Masi, A., Esposito, R., and Lebowitz, J.~L.}
\newblock Incompressible {N}avier-{S}tokes and {E}uler limits of the
  {B}oltzmann equation.
\newblock {\em Comm. Pure Appl. Math. 42}, 8 (1989), 1189--1214.

\bibitem{MR2569870}
{\sc Desvillettes, L., and Salvarani, F.}
\newblock Asymptotic behavior of degenerate linear transport equations.
\newblock {\em Bull. Sci. Math. 133}, 8 (2009), 848--858.

\bibitem{Desvillettes-Villani-2001}
{\sc Desvillettes, L., and Villani, C.}
\newblock On the trend to global equilibrium in spatially inhomogeneous
  entropy-dissipating systems: the linear {F}okker-{P}lanck equation.
\newblock {\em Comm. Pure Appl. Math. 54}, 1 (2001), 1--42.

\bibitem{DV02}
{\sc Desvillettes, L., and Villani, C.}
\newblock On a variant of {K}orn's inequality arising in statistical mechanics.
\newblock {\em ESAIM Control Optim. Calc. Var. 8\/} (2002), 603--619
  (electronic).
\newblock A tribute to J. L. Lions.

\bibitem{Desvillettes-Villani-2005}
{\sc Desvillettes, L., and Villani, C.}
\newblock On the trend to global equilibrium for spatially inhomogeneous
  kinetic systems: the {B}oltzmann equation.
\newblock {\em Invent. Math. 159}, 2 (2005), 245--316.

\bibitem{DiPernaLionsAnnals}
{\sc DiPerna, R.~J., and Lions, P.-L.}
\newblock On the {C}auchy problem for {B}oltzmann equations: global existence
  and weak stability.
\newblock {\em Ann. of Math. (2) 130}, 2 (1989), 321--366.

\bibitem{MR1088276}
{\sc DiPerna, R.~J., and Lions, P.-L.}
\newblock Global solutions of {B}oltzmann's equation and the entropy
  inequality.
\newblock {\em Arch. Rational Mech. Anal. 114}, 1 (1991), 47--55.

\bibitem{Dolbeault2009511}
{\sc Dolbeault, J., Mouhot, C., and Schmeiser, C.}
\newblock Hypocoercivity for kinetic equations with linear relaxation terms.
\newblock {\em Comptes Rendus Mathematique 347}, 9-10 (2009), 511 -- 516.

\bibitem{MR3324910}
{\sc Dolbeault, J., Mouhot, C., and Schmeiser, C.}
\newblock Hypocoercivity for linear kinetic equations conserving mass.
\newblock {\em Trans. Amer. Math. Soc. 367}, 6 (2015), 3807--3828.

\bibitem{MR2813582}
{\sc Duan, R.}
\newblock Hypocoercivity of linear degenerately dissipative kinetic equations.
\newblock {\em Nonlinearity 24}, 8 (2011), 2165--2189.

\bibitem{MR2966364}
{\sc Duan, R., and Li, W.-X.}
\newblock Hypocoercivity for the linear {B}oltzmann equation with confining
  forces.
\newblock {\em J. Stat. Phys. 148}, 2 (2012), 306--324.

\bibitem{DuanLiuSakamotoStrain}
{\sc Duan, R., Liu, S., Sakamoto, S., and Strain, R.~M.}
\newblock Global mild solutions of the {L}andau and non-cutoff {B}oltzmann
  equations.
\newblock Preprint arXiv:1904.12086.

\bibitem{MR0521262}
{\sc Duvaut, G., and Lions, J.-L.}
\newblock {\em Inequalities in mechanics and physics}.
\newblock Springer-Verlag, Berlin-New York, 1976.
\newblock Translated from the French by C. W. John, Grundlehren der
  Mathematischen Wissenschaften, 219.

\bibitem{MR1969727}
{\sc Eckmann, J.-P., and Hairer, M.}
\newblock Spectral properties of hypoelliptic operators.
\newblock {\em Comm. Math. Phys. 235}, 2 (2003), 233--253.

\bibitem{MR22750}
{\sc Friedrichs, K.~O.}
\newblock On the boundary-value problems of the theory of elasticity and
  {K}orn's inequality.
\newblock {\em Ann. of Math. (2) 48\/} (1947), 441--471.

\bibitem{GT}
{\sc Gallagher, I., and Tristani, I.}
\newblock On the convergence of smooth solutions from {B}oltzmann to
  {N}avier--{S}tokes.
\newblock {\em Annales Henri Lebesgue 3\/} (2020), 561--614.

\bibitem{GSR}
{\sc Golse, F., and Saint-Raymond, L.}
\newblock The {N}avier-{S}tokes limit of the {B}oltzmann equation for bounded
  collision kernels.
\newblock {\em Invent. Math. 155}, 1 (2004), 81--161.

\bibitem{MR775683}
{\sc Grisvard, P.}
\newblock {\em Elliptic problems in nonsmooth domains}, vol.~24 of {\em
  Monographs and Studies in Mathematics}.
\newblock Pitman (Advanced Publishing Program), Boston, MA, 1985.

\bibitem{Guo-2002-I}
{\sc Guo, Y.}
\newblock The {V}lasov-{P}oisson-{B}oltzmann system near {M}axwellians.
\newblock {\em Comm. Pure Appl. Math. 55}, 9 (2002), 1104--1135.

\bibitem{MR2679358}
{\sc Guo, Y.}
\newblock Decay and continuity of the {B}oltzmann equation in bounded domains.
\newblock {\em Arch. Ration. Mech. Anal. 197}, 3 (2010), 713--809.

\bibitem{Guoetal-specular2}
{\sc Guo, Y., Hwang, H.~J., Jang, J.~W., and Ouyang, Z.}
\newblock {$L^2$ decay for the linearized Landau equation with the specular
  boundary condition}.
\newblock Preprint arXiv:2009.01391.

\bibitem{MR4076068}
{\sc Guo, Y., Hwang, H.~J., Jang, J.~W., and Ouyang, Z.}
\newblock The {L}andau equation with the specular reflection boundary
  condition.
\newblock {\em Arch. Ration. Mech. Anal. 236}, 3 (2020), 1389--1454.

\bibitem{HamdacheARMA1992}
{\sc Hamdache, K.}
\newblock Initial-boundary value problems for the {B}oltzmann equation: global
  existence of weak solutions.
\newblock {\em Arch. Rational Mech. Anal. 119}, 4 (1992), 309--353.

\bibitem{MR3479064}
{\sc Han-Kwan, D., and L\'{e}autaud, M.}
\newblock Geometric analysis of the linear {B}oltzmann equation {I}. {T}rend to
  equilibrium.
\newblock {\em Ann. PDE 1}, 1 (2015), Art. 3, 84.

\bibitem{HN05}
{\sc Helffer, B., and Nier, F.}
\newblock {\em Hypoelliptic estimates and spectral theory for {F}okker-{P}lanck
  operators and {W}itten {L}aplacians}, vol.~1862 of {\em Lecture Notes in
  Mathematics}.
\newblock Springer-Verlag, Berlin, 2005.

\bibitem{MR2215889}
{\sc H\'{e}rau, F.}
\newblock Hypocoercivity and exponential time decay for the linear
  inhomogeneous relaxation {B}oltzmann equation.
\newblock {\em Asymptot. Anal. 46}, 3-4 (2006), 349--359.

\bibitem{Herau-Nier}
{\sc H{\'e}rau, F., and Nier, F.}
\newblock Isotropic hypoellipticity and trend to equilibrium for the
  {F}okker-{P}lanck equation with a high-degree potential.
\newblock {\em Arch. Ration. Mech. Anal. 171}, 2 (2004), 151--218.

\bibitem{JM}
{\sc Jiang, N., and Masmoudi, N.}
\newblock Boundary layers and incompressible {N}avier-{S}tokes-{F}ourier limit
  of the {B}oltzmann equation in bounded domain {I}.
\newblock {\em Comm. Pure Appl. Math. 70}, 1 (2017), 90--171.

\bibitem{JXZ}
{\sc Jiang, N., Xu, C.-J., and Zhao, H.}
\newblock Incompressible {N}avier-{S}tokes-{F}ourier limit from the {B}oltzmann
  equation: classical solutions.
\newblock {\em Indiana Univ. Math. J. 67}, 5 (2018), 1817--1855.

\bibitem{KM}
{\sc Kavian, O., and Mischler, S.}
\newblock The {F}okker-{P}lanck equation with subcritical confinement force.
\newblock \emph{arXiv:1512.07005}.

\bibitem{MR3762275}
{\sc Kim, C., and Lee, D.}
\newblock The {B}oltzmann equation with specular boundary condition in convex
  domains.
\newblock {\em Comm. Pure Appl. Math. 71}, 3 (2018), 411--504.

\bibitem{MR3840911}
{\sc Kim, C., and Lee, D.}
\newblock Decay of the {B}oltzmann equation with the specular boundary
  condition in non-convex cylindrical domains.
\newblock {\em Arch. Ration. Mech. Anal. 230}, 1 (2018), 49--123.

\bibitem{MR3199988}
{\sc Kuo, H.-W., Liu, T.-P., and Tsai, L.-C.}
\newblock Equilibrating effects of boundary and collision in rarefied gases.
\newblock {\em Comm. Math. Phys. 328}, 2 (2014), 421--480.

\bibitem{MR3626005}
{\sc L\'{e}autaud, M., and Lerner, N.}
\newblock Energy decay for a locally undamped wave equation.
\newblock {\em Ann. Fac. Sci. Toulouse Math. (6) 26}, 1 (2017), 157--205.

\bibitem{Lebeau}
{\sc Lebeau, G.}
\newblock \'{E}quation des ondes amorties.
\newblock In {\em Algebraic and geometric methods in mathematical physics
  ({K}aciveli, 1993)}, vol.~19 of {\em Math. Phys. Stud.} Kluwer Acad. Publ.,
  Dordrecht, 1996, pp.~73--109.

\bibitem{LebeauRobb}
{\sc Lebeau, G., and Robbiano, L.}
\newblock Stabilisation de l'\'equation des ondes par le bord.
\newblock {\em Duke Math. J. 86}, 3 (1997), 465--491.

\bibitem{MR1295942}
{\sc Lions, P.-L.}
\newblock Compactness in {B}oltzmann's equation via {F}ourier integral
  operators and applications. {III}.
\newblock {\em J. Math. Kyoto Univ. 34}, 3 (1994), 539--584.

\bibitem{MSR}
{\sc Masmoudi, N., and Saint-Raymond, L.}
\newblock From the {B}oltzmann equation to the {S}tokes-{F}ourier system in a
  bounded domain.
\newblock {\em Comm. Pure Appl. Math. 56}, 9 (2003), 1263--1293.

\bibitem{MischlerCMP2000}
{\sc Mischler, S.}
\newblock On the initial boundary value problem for the
  {V}lasov-{P}oisson-{B}oltzmann system.
\newblock {\em Comm. Math. Phys. 210}, 2 (2000), 447--466.

\bibitem{MR1765137}
{\sc Mischler, S.}
\newblock On the trace problem for solutions of the {V}lasov equation.
\newblock {\em Comm. Partial Differential Equations 25}, 7-8 (2000),
  1415--1443.

\bibitem{MischlerENS2010}
{\sc Mischler, S.}
\newblock Kinetic equations with {M}axwell boundary conditions.
\newblock {\em Ann. Sci. \'{E}c. Norm. Sup\'{e}r. (4) 43}, 5 (2010), 719--760.

\bibitem{Mouhot-Neumann}
{\sc Mouhot, C., and Neumann, L.}
\newblock Quantitative perturbative study of convergence to equilibrium for
  collisional kinetic models in the torus.
\newblock {\em Nonlinearity 19}, 4 (2006), 969--998.

\bibitem{Rachid}
{\sc Rachid, M.}
\newblock Incompressible {N}avier-{S}tokes-{F}ourier limit from the {L}andau
  equation.
\newblock Preprint hal-03035871.

\bibitem{RockWang}
{\sc R{\"o}ckner, M., and Wang, F.-Y.}
\newblock Weak {P}oincar\'e inequalities and {$L^2$}-convergence rates of
  {M}arkov semigroups.
\newblock {\em J. Funct. Anal. 185}, 2 (2001), 564--603.

\bibitem{SR-book}
{\sc Saint-Raymond, L.}
\newblock {\em Hydrodynamic limits of the {B}oltzmann equation}, vol.~1971 of
  {\em Lecture Notes in Mathematics}.
\newblock Springer-Verlag, Berlin, 2009.

\bibitem{MR1452171}
{\sc Savar\'{e}, G.}
\newblock Regularity and perturbation results for mixed second order elliptic
  problems.
\newblock {\em Comm. Partial Differential Equations 22}, 5-6 (1997), 869--899.

\bibitem{GS1}
{\sc Strain, R.~M., and Guo, Y.}
\newblock Almost exponential decay near {M}axwellian.
\newblock {\em Comm. Partial Differential Equations 31}, 1-3 (2006), 417--429.

\bibitem{GS2}
{\sc Strain, R.~M., and Guo, Y.}
\newblock Exponential decay for soft potentials near {M}axwellian.
\newblock {\em Arch. Ration. Mech. Anal. 187}, 2 (2008), 287--339.

\bibitem{Mem-villani}
{\sc Villani, C.}
\newblock Hypocoercivity.
\newblock {\em Mem. Amer. Math. Soc. 202}, 950 (2009), iv+141.

\end{thebibliography}

\bigskip\bigskip

\signab

\signkc

\signsm

\signit

\end{document}